\documentclass{amsart}
\usepackage{amsmath,amssymb,amsfonts,epsfig}

\newtheorem{theorem}{Theorem}[section]
\newtheorem{proposition}[theorem]{Proposition}

\newtheorem{conjecture}[theorem]{Conjecture}

\theoremstyle{definition}

\newtheorem{problem}[theorem]{Problem}

\theoremstyle{remark}
\newtheorem{remark}[theorem]{Remark}

\numberwithin{equation}{section}

\def\DJ{{\hbox{D\kern-.8em\raise.15ex\hbox{--}\kern.35em}}}
\def\DJo{\DJ okovi\'c}

\providecommand{\bysame}{\leavevmode\hbox to3em{\hrulefill}\thinspace}

\def\vf{{\varphi}}
\def\bC{{\bf C}}

\def\bZ{{\bf Z}}
\def\pA{{\mathcal A}}

\def\pP{{\mathcal P}}

\def\tr{{\rm tr\;}}

\def\GL{{\mbox{\rm GL}}}
\def\SL{{\mbox{\rm SL}}}

\begin{document}

\title[Poincar\'{e} series of some trace algebras]
{Poincar\'{e} series of some pure and mixed trace algebras
of two generic matrices}

\author[D.\v{Z}. \DJo{}]
{Dragomir \v{Z}. \DJo{}}

\address{Department of Pure Mathematics, University of Waterloo,
Waterloo, Ontario, N2L 3G1, Canada}

\email{djokovic@uwaterloo.ca}


\thanks{The author was supported in part by the NSERC
Grant A-5285.}

\date{}

\begin{abstract}
We work over a field $K$ of characteristic zero.
The Poincar\'{e} series for the algebra $C_{n,2}$ of $\GL_n$-invariants
and the algebra $T_{n,2}$ of $\GL_n$-concomitants of
two generic $n\times n$ matrices $x$ and $y$ are presented for $n\le6$.
Both simply graded and bigraded cases are included.
The cases $n\le4$ were known previously. If $n=5$ or 6, we show
that $C_{n,2}$ has no bigraded system of parameters.

For the algebra $C_{4,2}$ we construct a minimal set of generators,
each of the form $\tr(w(x,y))$ where $w$ is a word in two letters,
and give an application to Specht's theorem on unitary similarity
of two complex matrices. We also construct a minimal set of generators
for the algebra $C_{5,2}$.
It consists of 173 bihomogeneous polynomials. Five of them are
symmetric (i.e., satisfy $f(x,y)=f(y,x)$), four are skew-symmetric
and the remaining 164 generators occur in pairs
$\{f(x,y),f(y,x)\}$.

By identifying the space $M_n^2$ of pairs of $n\times n$ matrices
with $M_n\otimes K^2$, we extend the action of 
$\GL_n$ to $\GL_n\times\GL_2$.
For $n\le5$, we compute the Poincar\'{e} series for the
polynomial invariants of this action when restricted to the subgroups
$\GL_n\times\SL_2$ and $\GL_n\times\Delta_1$, where $\Delta_1$
is the maximal torus of $\SL_2$ consisting of diagonal matrices.

Five conjectures are proposed concerning the lowest term numerators
and denominators of various Poincar\'{e} series mentioned above.
Some heuristic formulas and open problems are indicated.
\end{abstract}

\maketitle

\section{Introduction}

Let $K$ be a field of characteristic 0, $M_n=M_n(K)$ the
$K$-algebra of $n\times n$ matrices over $K$, and $d$ a positive integer.
The general linear group $\GL_n=\GL_n(K)$ acts on the direct product
$M_n^d$ of $d$ copies of $M_n$ by simultaneous conjugation
\[ a\cdot(x_1,\ldots,x_d)=(ax_1a^{-1},\ldots,ax_d a^{-1}). \]
This gives rise to an action of $\GL_n$ on the algebra 
$K[M_n^d]$ of polynomial functions on $M_n^d$.
We shall view the entries of $x_k$, $k=1,\ldots,d$, also as linear functions
on $M_n^d$, thus providing a convenient system of coordinates on this space.

It is a well known fact, due to Procesi \cite{CPr} and Razmyslov \cite{YR},
that the algebra $C_{n,d}=K[M_n^d]^{\GL_n}$ of $\GL_n$-invariants in
$K[M_n^d]$ is generated by all traces
\[ \tr(z_1z_2\cdots z_k),\quad z_1,z_2,\ldots,z_k\in\{x_1,\ldots,x_d\},
\quad k\ge1. \]
This assertion remains valid if one imposes the restriction
$k\le n^2$. For these facts and many other known properties
of $C_{n,d}$ we refer the interested reader to one of the references
\cite{DF,Dr1,EF3} and the papers quoted there.

Similarly, $\GL_n$ acts on the (noncommutative) algebra of polynomial
maps $M_n^d\to M_n$. Its subalgebra consisting of $\GL_n$-equivariant
maps will be denoted by $T_{n,d}$.
The algebra $C_{n,d}$ resp. $T_{n,d}$ is known as the {\em pure}
resp. {\em mixed trace algebra} of $d$ generic $n\times n$ matrices.

By assigning the degree $(1,0,\ldots,0)$ to the entries of the matrix $x_1$,
the degree $(0,1,0,\ldots,0)$ to the entries of $x_2$, etc., one obtains a
$\bZ^d$-gradation of the algebras $C_{n,d}$ and $T_{n,d}$.
The total degree provides these algebras with the ordinary
$\bZ$-gradation. In this paper we are mostly interested in the case
$d=2$. Until further notice, we assume that this is the case and
we set $x=x_1$ and $y=x_2$.

We shall now describe briefly the contents of each section.
The simply graded and bigraded Poincar\'{e}
series for both $C_{n,2}$ and $T_{n,2}$ have been explicitly
computed for $n\le4$ (see \cite{BS,T1,T2,T3}).

In Section \ref{z2-grad} resp. \ref{z1-grad} we present the results of
our computations giving the bigraded resp. simply graded Poincar\'{e} 
series for the algebras $C_{n,2}$ and $T_{n,2}$ for $n=5,6$.
For the bigraded case see Theorem \ref{Prva} and Tables 1 and 2.
For the simply graded case see Tables 3 and 4, as well as Tables 9
and 11 in Appendices C and D, respectively.
We also show ( see Proposition \ref{nema-Z2} )
that the $\bZ^2$-graded algebra $C_{n,2}$ has no
bigraded system of parameters if $n=5,6$.

In Section \ref{gen-c42} we construct a minimal set of generators
of $C_{4,2}$ (see Theorem \ref{Gen-1})
which consists of 32 elements of the form
$\tr(w(x,y))$, where $w(x,y)$ is a word in the matrices $x$ and $y$.
This problem has been already solved by Drensky and Sadikova
\cite{DS}, but their generators are not of this simple form.
However, they claim that their choice is better suited for
finding a presentation of $C_{4,2}$, which is apparently still
an open problem.

We give an application to the problem of unitary similarity of
two complex $4\times 4$ matrices. Recall that Pearcy
\cite{CP} showed more than 40 years ago that two $3\times3$
complex matrices $a$ and $b$ are unitarily similar iff the
equalities $\tr(w(a,a^*))=\tr(w(b,b^*))$ hold when $w$ runs 
through his list of nine words. Sibirski\u\i\ \cite{KS} showed later
that two of these nine words are redundant. We
provide a similar test for unitary similarity of two
$4\times 4$ complex matrices $a$ and $b$.
Some of the words in our generating set occur in pairs
$\{w_1,w_2\}$ such that $\tr(w_1(a,a^*))$ and $\tr(w_2(a,a^*))$ are
complex conjugates of each other for any complex matrix $a$.
Hence we can omit a word from any such pair. In this way
we obtain a reduced list of only 20 words which suffices for
testing the unitary similarity of $4\times4$ complex matrices
(see Theorem \ref{Un-sim}).

In Section \ref{gen-c52} we construct a minimal set of generators,
$\pP$, for the algebra $C_{5,2}$.
It consists of 173 bihomogeneous polynomials.

In Section \ref{veca-gr} we identify the space $M_n^d$ with the
tensor product $M_n\otimes K^d$ and extend the action of $\GL_n$
to $\GL_n\times\GL_d$ by letting $\GL_d$ act on $K^d$ by
multiplication. We denote by
$C_{n,d}^\#$ the subalgebra of $K[M_n^d]$ consisting of
$\GL_n\times\SL_d$-invariant functions.
In Table 7 we record the results of our computation of
the Poincar\'{e} series of $C_{n,2}^\#$, $n\le5$.

In Section \ref{manja-gr} we restrict the action of
$\GL_n\times\SL_d$ to $\GL_n\times\Delta_{d-1}$,
where $\Delta_{d-1}$ is the maximal torus of $\SL_d$
consisting of diagonal matrices.
We denote by $C_{n,d}^\bullet$ the subalgebra of $K[M_n^d]$
consisting of $\GL_n\times\Delta_{d-1}$-invariant functions.
In Table 8 we record the Poincar\'{e} series for the algebras
$C_{n,2}^\bullet$, $n\le5$.

In Section \ref{hipot} we propose four conjectures 
concerning the lowest term numerators and denominators
of the Poincar\'{e} series of $C_{n,2}$ and $T_{n,2}$,
one more conjecture about the series for
$C_{n,2}^\#$ and $C_{n,2}^\bullet$, and state an open problem.

The numerators of the bigraded Poincar\'{e} series of $C_{6,2}$ 
and $T_{6,2}$ have 1169 and 854 terms, respectively. For that
reason we give them separately in Appendix A. By using the fact
that these polynomials are symmetric in the two variables and
that they satisfy a simple functional equation, we can write
them in a compact form.

Appendix B gives some details concerning the verification that our
expression for $P(C_{5,2};s,t)$ agrees 
with Formanek's expansion in terms of Schur functions (see \cite{BS}).

Appendix C contains the Taylor expansions
(listing all coefficients of degree less than 20)
of $P(C_{n,2};t)$ for $n\le6$ and for the closely
related algebras $C_{n,2}(0)$, introduced in Section \ref{gen-c52}.
We also tabulate some of these coefficients, make a couple of
interesting observations, and indulge in some speculative thinking.

Appendix D treats in the same way the Taylor expansions of
$P(T_{n,2};t)$ for $n\le6$.

We warn the reader that several assertions and formulae that appear
in the last two appendices are of hypothetical character and are
included there only as a suggestion to the interested reader as
deserving further consideration and study.

I thank my student K.-C. Chan for his comments on several preliminary
versions of this paper.

\section{ Bigraded Poincar\'{e} series }
\label{z2-grad}

Let $P(C_{n,2};s,t)$ resp. $P(T_{n,2};s,t)$
denote the bigraded Poincar\'{e} series of the pure resp. mixed
trace algebra $C_{n,2}$ resp. $T_{n,2}$.
It is known that these series are given by symmetric rational
functions in the two variables $s$ and $t$.

We can write these rational functions in lowest terms as
\begin{equation} \label{poen_Cst}
P(C_{n,2};s,t) = \frac{ N(C_{n,2};s,t)}{ D(C_{n,2};s,t)}
\end{equation}
and
\begin{equation} \label{poen_Tst}
P(T_{n,2};s,t) = \frac{ N(T_{n,2};s,t)}{ D(T_{n,2};s,t)}, 
\end{equation}
where we normalize $N$ and $D$ by demanding that they both
have constant term 1.

We recall that Van den Bergh \cite[Proposition 5.1]{VB} has shown 
(in a more general setting)
that these two rational functions can be written as fractions
whose numerator is a polynomial with integer coefficients
and the denominator is a product
of terms $1-u$ with $u$ a monomial in $s$ and $t$.
He does not claim that the latter property holds if the fraction is
in lowest terms. However, this is indeed true for $n\le6$.

The Poincar\'{e} series (\ref{poen_Cst}) and (\ref{poen_Tst})
have been computed by Teranishi \cite{T1,T2,T3} for $n\le4$.
The formulae for the case $n=4$ were computed independently
by Berele and Stembridge \cite{BS}. They have
also corrected several misprints in Teranishi's formulas in \cite{T2}.

For the sake of completeness and the reader's convenience, we
give explicit formulae for
the numerators and denominators for all $n\le6$.
For $n\le5$, the numerators
$N(C_{n,2};s,t)$ resp. $N(T_{n,2};s,t)$
are given in Table 1 resp. 2. For $n=6$,
the numerators are huge and are given in Appendix A.

We remark that these numerators are symmetric polynomials
in $s$ and $t$, and that they also satisfy the functional equation
\[ (st)^d N(C_{n,2};s^{-1},t^{-1})=N(C_{n,2};s,t), \]
where $d$ is the degree of
$N(C_{n,2};s,t)$ as a polynomial in $s$.
For that reason (and to save space) there is no need to write all
the terms of $N(C_{n,2};s,t)$.
This remark also applies to the numerators $N(T_{n,2};s,t)$.

\begin{center}
{\bf Table 1: Numerators $N(C_{n,2};s,t)$}

\begin{eqnarray*}
&& N(C_{1,2};s,t)=N(C_{2,2};s,t)=1, \\
&& N(C_{3,2};s,t)= 1-st+{s}^{2}{t}^{2}\, , \\
&& N(C_{4,2};s,t)=\\
&&  \left( 1-st+{s}^{2}{t}^{2} \right)
 \left( 1-s{t}^{2}-{s}^{2}t+{s}^{2}{t}^{3}+{s}^{3}{t}^{2}
+{s}^{2}{t}^{4}+2\,{s}^{3}{t}^{3}+{s}^{4}{t}^{2}+{s}^{3}{t}^{4} \right. \\
&& \left. +{s}^{4}{t}^{3}-{s}^{4}{t}^{5}-{s}^{5}{t}^{4}
+{s}^{6}{t}^{6} \right) \\
&& =1-st-s{t}^{2}-{s}^{2}t+{s}^{2}{t}^{2}+2\,{s}^{2}{t}^{3}
+2\,{s}^{3}{t}^{2}+{s}^{2}{t}^{4}+2\,{s}^{3}{t}^{3}+{s}^{4}{t}^{2}
-{s}^{3}{t}^{4} \\
&& -{s}^{4}{t}^{3}-{s}^{3}{t}^{5}-2\,{s}^{4}{t}^{4}
-{s}^{5}{t}^{3}- \cdots -{s}^{7}{t}^{7}+{s}^{8}{t}^{8}, \\
&& N(C_{5,2};s,t)= \\
&& 1-st-2\,s{t}^{2}-2\,{s}^{2}t-s{t}^{3}-{s}^{3}t+3\,{s}^{2}{t}^{3}
+3\,{s}^{3}{t}^{2}+4\,{s}^{2}{t}^{4}+7\,{s}^{3}{t}^{3}+4\,{s}^{4}{t}^{2} \\
&& +3\,{s}^{2}{t}^{5}+4\,{s}^{3}{t}^{4}+4\,{s}^{4}{t}^{3}+3\,{s}^{5}{t}^{2}
+{s}^{2}{t}^{6}-2\,{s}^{3}{t}^{5}-5\,{s}^{4}{t}^{4}-2\,{s}^{5}{t}^{3}+{s}^{6}{t}^{2} \\
&& -5\,{s}^{3}{t}^{6}-13\,{s}^{4}{t}^{5}-13\,{s}^{5}{t}^{4}-5\,{s}^{6}{t}^{3}
-4\,{s}^{3}{t}^{7}-10\,{s}^{4}{t}^{6}-11\,{s}^{5}{t}^{5}-10\,{s}^{6}{t}^{4} \\
&& -4\,{s}^{7}{t}^{3}-2\,{s}^{3}{t}^{8}-2\,{s}^{4}{t}^{7}+3\,{s}^{5}{t}^{6}
+3\,{s}^{6}{t}^{5}
-2\,{s}^{7}{t}^{4}-2\,{s}^{8}{t}^{3}+5\,{s}^{4}{t}^{8}+18\,{s}^{5}{t}^{7} \\
&& +27\,{s}^{6}{t}^{6}+18\,{s}^{7}{t}^{5}+5\,{s}^{8}{t}^{4}+4\,{s}^{4}{t}^{9}
+18\,{s}^{5}{t}^{8}+29\,{s}^{6}{t}^{7}+29\,{s}^{7}{t}^{6}+18\,{s}^{8}{t}^{5} \\
&& +4\,{s}^{9}{t}^{4}+2\,{s}^{4}{t}^{10}+8\,{s}^{5}{t}^{9}+12\,{s}^{6}{t}^{8}
+11\,{s}^{7}{t}^{7}+12\,{s}^{8}{t}^{6}+8\,{s}^{9}{t}^{5}+2\,{s}^{10}{t}^{4} \\
&& -{s}^{5}{t}^{10}-11\,{s}^{6}{t}^{9}-24\,{s}^{7}{t}^{8}-24\,{s}^{8}{t}^{7}
-11\,{s}^{9}{t}^{6}-{s}^{10}{t}^{5}-3\,{s}^{5}{t}^{11}-17\,{s}^{6}{t}^{10} \\
&& -40\,{s}^{7}{t}^{9}-47\,{s}^{8}{t}^{8}-40\,{s}^{9}{t}^{7}-17\,{s}^{10}{t}^{6}
-3\,{s}^{11}{t}^{5}-{s}^{5}{t}^{12}-10\,{s}^{6}{t}^{11}-24\,{s}^{7}{t}^{10} \\
&& -33\,{s}^{8}{t}^{9}-33\,{s}^{9}{t}^{8}-24\,{s}^{10}{t}^{7}-10\,{s}^{11}{t}^{6}
-{s}^{12}{t}^{5}-{s}^{6}{t}^{12}-{s}^{7}{t}^{11}+9\,{s}^{8}{t}^{10} \\
&& +13\,{s}^{9}{t}^{9}+9\,{s}^{10}{t}^{8}-{s}^{11}{t}^{7}-{s}^{12}{t}^{6}
+{s}^{6}{t}^{13}+10\,{s}^{7}{t}^{12}+34\,{s}^{8}{t}^{11}+54\,{s}^{9}{t}^{10} \\
&& +54\,{s}^{10}{t}^{9}+34\,{s}^{11}{t}^{8}+10\,{s}^{12}{t}^{7}+{s}^{13}{t}^{6}
+6\,{s}^{7}{t}^{13}+27\,{s}^{8}{t}^{12}+49\,{s}^{9}{t}^{11} \\
&& +62\,{s}^{10}{t}^{10}+49\,{s}^{11}{t}^{9}+27\,{s}^{12}{t}^{8}+6\,{s}^{13}{t}^{7}
+5\,{s}^{8}{t}^{13}+10\,{s}^{9}{t}^{12}+13\,{s}^{10}{t}^{11} \\
&& +13\,{s}^{11}{t}^{10}+10\,{s}^{12}{t}^{9}+5\,{s}^{13}{t}^{8}-4\,{s}^{8}{t}^{14}
-19\,{s}^{9}{t}^{13}-37\,{s}^{10}{t}^{12}-45\,{s}^{11}{t}^{11} \\
&& -37\,{s}^{12}{t}^{10}-19\,{s}^{13}{t}^{9}-4\,{s}^{14}{t}^{8}-2\,{s}^{8}{t}^{15}
-20\,{s}^{9}{t}^{14}-48\,{s}^{10}{t}^{13}-72\,{s}^{11}{t}^{12} \\
&& -72\,{s}^{12}{t}^{11}-48\,{s}^{13}{t}^{10}-20\,{s}^{14}{t}^{9}
-2\,{s}^{15}{t}^{8}-\,\cdots\,-{s}^{22}{t}^{22}+{s}^{23}{t}^{23}. 
\end{eqnarray*}

\end{center}

One can easily verify that all coefficients of the extended numerators
\begin{eqnarray*}
&& N^*(C_{3,2};s,t)= \left( 1+st \right) N(C_{3,2};s,t), \\
&& N^*(C_{4,2};s,t)= \left( 1+st \right) \left( 1+s{t}^{2} \right)
\left( 1+{s}^{2}t \right) N(C_{4,2};s,t)
\end{eqnarray*}
are nonnegative integers. The corresponding extended denominators
\begin{eqnarray*}
&& D^*(C_{3,2};s,t)= \left( 1+st \right) D(C_{3,2};s,t), \\
&& D^*(C_{4,2};s,t)= \left( 1+st \right) \left( 1+s{t}^{2} \right)
\left( 1+{s}^{2}t \right) D(C_{4,2};s,t)
\end{eqnarray*}
can still be written as products of the binomials of the form $1-u$.
Explicitly, we have
\begin{eqnarray*}
N^*(C_{3,2};s,t) &=& 1+{s}^{3}{t}^{3}, \\
D^*(C_{3,2};s,t) &=& \left( 1-s \right) \left( 1-t \right) \cdot 
 \left( 1-{s}^{2} \right) \left( 1-st \right)
\left( 1-{t}^{2} \right) \cdot \\
&& \left( 1-{s}^{3} \right)  \left( 1-{s}^{2}t \right)
	\left( 1-s{t}^{2} \right) \left( 1-{t}^{3} \right) \cdot 
\left( 1-{s}^{2}{t}^{2} \right), \\
N^*(C_{4,2};s,t) &=&
1+{s}^{2}{t}^{3}+{s}^{3}{t}^{2}+2\,{s}^{3}{t}^{3}+{s}^{3}{t}^{4}+{s}^{
4}{t}^{3}+{s}^{3}{t}^{5}+2\,{s}^{4}{t}^{4}+{s}^{5}{t}^{3} \\
&& +{s}^{3}{t}^{6}+{s}^{4}{t}^{5}+{s}^{5}{t}^{4}+{s}^{6}{t}^{3}
+{s}^{4}{t}^{6}+2\,{s}^{5}{t}^{5}+{s}^{6}{t}^{4}+2\,{s}^{5}{t}^{6} \\
&& +2\,{s}^{6}{t}^{5}+2\,{s}^{6}{t}^{6}+ \cdots 
+{s}^{9}{t}^{10}+{s}^{10}{t}^{9}+{s}^{12}{t}^{12}, \\
D^*(C_{4,2};s,t) &=& \left( 1-s \right) \left( 1-t \right) \cdot 
\left( 1-{s}^{2} \right) \left( 1-st \right) \left( 1-{t}^{2} \right) \cdot \\
&& \left( 1-{s}^{3} \right)  \left( 1-{s}^{2}t \right)
	\left( 1-s{t}^{2} \right) \left( 1-{t}^{3} \right) \cdot \\
&& \left( 1-{s}^{4} \right) \left( 1-{s}^{3}t \right) 
	\left( 1-{s}^{2}{t}^{2} \right)^{2} \left( 1-s{t}^{3} \right)
	\left( 1-{t}^{4} \right) \cdot \\
&& \left( 1-{s}^{4}{t}^{2} \right) \left( 1-{s}^{2}{t}^{4} \right).
\end{eqnarray*}

On the other hand, the Poincar\'{e} series
$P(C_{5,2};s,t)$ and $P(C_{6,2};s,t)$
do not admit similar expressions because
$N(C_{5,2};1,1)=N(C_{6,2};1,1)=0$.

\begin{proposition} \label{nema-Z2}
The algebras $C_{5,2}$ and $C_{6,2}$ have no bigraded system of parameters.
\end{proposition}

\begin{proof}
Assume that $C_{5,2}$ has a bigraded system of parameters, say
$\{P_1,\ldots,P_{26}\}$. Since $C_{5,2}$ is a Cohen--Macauley algebra,
it is a free graded module over the polynomial algebra
$K[P_1,\ldots,P_{26}]$. Consequently, this module has a free
bigraded basis, say $\{Q_1=1,\; Q_2,\ldots,Q_m\}$. It follows that
\[ P(C_{5,2};s,t)=\frac
{\sum_{j=1}^m s^{e'_j} t^{e''_j}}
{\prod_{i=1}^{26} \left( 1-s^{d'_i} t^{d''_i} \right) }=
\frac{N(C_{5,2};s,t)}{D(C_{5,2};s,t)}, \]
where $(d'_i,d''_i)$ resp. $(e'_j,e''_j)$ is the bidegree of
$P_i$ resp $Q_j$. As
$N(C_{5,2};s,t)$ and $D(C_{5,2};s,t)$ are relatively prime,
we deduce that
\[ \sum_{j=1}^m s^{e'_j} t^{e''_j}=
N(C_{5,2};s,t)R(s,t), \]
where $R(s,t)$ is a polynomial in $s$ and $t$.
By setting $s=t=1$ in this identity, we obtain the contradiction $m=0$
since $N(C_{5,2};1,1)=0$.

The proof for the algebra $C_{6,2}$ is similar.
\end{proof}

\newpage
\begin{center}
{\bf Table 2: Numerators $N(T_{n,2};s,t)$}

\begin{eqnarray*}
 && N(T_{1,2};s,t)=N(T_{2,2};s,t)=N(T_{3,2};s,t)=1, \\
 && N(T_{4,2};s,t)=1+{s}^{2}{t}^{2}+{s}^{2}{t}^{3}+{s}^{3}{t}^{2}
+{s}^{3}{t}^{3}+{s}^{5}{t}^{5}, \\ 
 && N(T_{5,2};s,t)= \\
&& 1-s{t}^{2}-{s}^{2}t+{s}^{2}{t}^{2}+2\,{s}^{2}{t}^{3}+2\,{s}^{3}{t}^{2}
+2\,{s}^{2}{t}^{4}+5\,{s}^{3}{t}^{3}+2\,{s}^{4}{t}^{2}+{s}^{2}{t}^{5}
+2\,{s}^{3}{t}^{4} \\
&& +2\,{s}^{4}{t}^{3}+{s}^{5}{t}^{2}-{s}^{3}{t}^{6}-3\,{s}^{4}{t}^{5}
-3\,{s}^{5}{t}^{4}-{s}^{6}{t}^{3}-{s}^{3}{t}^{7}-{s}^{4}{t}^{6}+{s}^{5}{t}^{5}
-{s}^{6}{t}^{4}-{s}^{7}{t}^{3} \\
&& +4\,{s}^{5}{t}^{6}+4\,{s}^{6}{t}^{5}+{s}^{4}{t}^{8}+6\,{s}^{5}{t}^{7}
+12\,{s}^{6}{t}^{6}+6\,{s}^{7}{t}^{5}+{s}^{8}{t}^{4}+3\,{s}^{5}{t}^{8}
+6\,{s}^{6}{t}^{7} \\
&& +6\,{s}^{7}{t}^{6}+3\,{s}^{8}{t}^{5}+{s}^{6}{t}^{8}+{s}^{7}{t}^{7}+{s}^{8}{t}^{6}
-{s}^{5}{t}^{10}-5\,{s}^{6}{t}^{9}-9\,{s}^{7}{t}^{8}-9\,{s}^{8}{t}^{7}
-5\,{s}^{9}{t}^{6} \\
&& -{s}^{10}{t}^{5}-4\,{s}^{6}{t}^{10}-9\,{s}^{7}{t}^{9}-7\,{s}^{8}{t}^{8}
-9\,{s}^{9}{t}^{7}-4\,{s}^{10}{t}^{6}-{s}^{6}{t}^{11}
-5\,{s}^{7}{t}^{10}-5\,{s}^{8}{t}^{9} \\
&& -5\,{s}^{9}{t}^{8}-5\,{s}^{10}{t}^{7}-{s}^{11}{t}^{6}+3\,{s}^{8}{t}^{10}
+6\,{s}^{9}{t}^{9}+3\,{s}^{10}{t}^{8}+4\,{s}^{8}{t}^{11}+6\,{s}^{9}{t}^{10} \\
&& +6\,{s}^{10}{t}^{9}+4\,{s}^{11}{t}^{8}+\,\cdots\,-{s}^{17}{t}^{18}
-{s}^{18}{t}^{17}+{s}^{19}{t}^{19}. 
\end{eqnarray*}

\end{center}

The denominators
$D(C_{n,2};s,t)$ and $D(T_{n,2};s,t)$
are closely related to the product
\begin{equation} \label{Pi}
\Pi_n(s,t) = \prod_{i=1}^n (1-s^i)(1-t^i)
\prod_{j=1}^{i-1} (1-s^{i-j}t^j)^{\min(i,n+1-i)}.
\end{equation}

\begin{theorem} \label{Prva}
For $n\le5$, $N(C_{n,2};s,t)$ and $N(T_{n,2};s,t)$
are given by Tables 1 and 2, respectively.
For $n=6$ they are given in Appendix A.
For $n\le6$, $D(C_{n,2};s,t)$ and $D(T_{n,2};s,t)$ are given by
\begin{eqnarray*}
&& D(C_{n,2};s,t)=\Pi_n(s,t),\quad n\le5; \\
&& D(C_{6,2};s,t)=(1-st)\Pi_6(s,t); \\
&& D(T_{n,2};s,t)=(1+s+\cdots+s^{n-1})^{-1}
(1+t+\cdots+t^{n-1})^{-1} D(C{_n,2};s,t).
\end{eqnarray*}
\end{theorem}

For instance, we have
\begin{eqnarray*}
D(C_{5,2};s,t) &=& \left( 1-s \right) \left( 1-t \right) \cdot
\left( 1-{s}^{2} \right) \left( 1-st \right)^{2}
\left( 1-{t}^{2} \right) \cdot \\
&& \left( 1-{s}^{3} \right)  \left( 1-{s}^{2}t \right)^{3}
	\left( 1-s{t}^{2} \right)^{3} \left( 1-{t}^{3} \right) \cdot \\
&& \left( 1-{s}^{4} \right) \left( 1-{s}^{3}t \right)^{2} 
	\left( 1-{s}^{2}{t}^{2} \right)^{2} \left( 1-s{t}^{3} \right)^{2}  
	\left( 1-{t}^{4} \right) \cdot \\
&& \left( 1-{s}^{5} \right)  \left( 1-{s}^{4}t \right) 
	\left( 1-{s}^{3}{t}^{2} \right) \left( 1-{s}^{2}{t}^{3} \right) 
	\left( 1-s{t}^{4} \right) \left( 1-{t}^{5}\right), \\
D(T_{5,2};s,t) &=& \left( 1-s \right)^{2} \left( 1-t \right)^{2} \cdot 
\left( 1-{s}^{2} \right) \left( 1-st \right) ^{2}
\left( 1-{t}^{2} \right) \cdot \\
&& \left( 1-{s}^{3} \right)  \left( 1-{s}^{2}t \right)^{3}
	\left( 1-s{t}^{2} \right)^{3} \left( 1-{t}^{3} \right) \cdot \\
&& \left( 1-{s}^{4} \right) \left( 1-{s}^{3}t \right)^{2}
	\left( 1-{s}^{2}{t}^{2} \right)^{2} \left( 1-s{t}^{3} \right)^{2}
	\left( 1-{t}^{4} \right) \cdot \\
&&  \left( 1-{s}^{4}t \right)  \left( 1-{s}^{3}{t}^{2} \right)
 \left( 1-{s}^{2}{t}^{3} \right) \left( 1-s{t}^{4} \right).
\end{eqnarray*}

\begin{proof}
By using the well known Molien--Weyl formula (see \cite{DK}), we have:
\begin{eqnarray*}
P(C_{n,2};s,t) &=& \frac{1}{(1-s)^n(1-t)^n} \cdot
\frac{1}{(2\pi i)^{n-1}} \cdot \\
&& \int_{|x_1|=1}\cdots\int_{|x_{n-1}|=1}
\prod_{1\le k\le r\le n-1}
\frac{1-x_kx_{k+1}\cdots x_r}{\vf_{k,r}}
\frac{{\rm d}x_{n-1}}{x_{n-1}}\cdots
\frac{{\rm d}x_1}{x_1},
\end{eqnarray*}
where
\begin{eqnarray*}
\vf_{k,r} &=& (1-sx_k x_{k+1}\cdots x_r)
(1-tx_k x_{k+1}\cdots x_r)\cdot \\
&& (1-s(x_k x_{k+1}\cdots x_r)^{-1})
(1-t(x_k x_{k+1}\cdots x_r)^{-1}),
\end{eqnarray*}
the integration is performed over the unit circles
(in the counterclockwise direction),
and the variables $s$ and $t$ have small moduli.

A similar formula is valid for $P(T_{n,2};s,t)$.
One has just to multiply the above integrand by the function
\[ n+\sum_{r=1}^{n-1} \sum_{k=1}^r \left(
x_k x_{k+1}\cdots x_r+\frac{1}{x_k x_{k+1}\cdots x_r} \right). \]

For $n=5$ and $n=6$ we have computed the two types of integrals
by using MAPLE \cite{Map}.
Each of the cases $n=6$ required about two weeks of computing time
on a machine running R10000 CPU at 250 MHz with 8 GB of RAM.
\end{proof}

We have verified independently the low degree ($\le25$) coefficients
in the Taylor expansions of $P(C_{5,2};s,t)$ and $P(C_{6,2};s,t)$ 
by using a formula due to Formanek (see Appendix B).

\section{ Simply graded Poincar\'{e} series }
\label{z1-grad}

Let $P(C_{n,2};t)$ denote the simply graded Poincar\'{e} series of
$C_{n,2}$ and $P(T_{n,2};t)$ the one for $T_{n,2}$.
Of course, one has
\begin{equation} \label{veze}
P(C_{n,2};t)=P(C_{n,2};t,t),\quad P(T_{n,2};t)=P(T_{n,2};t,t).
\end{equation}
By setting $s=t$ in the integral formula for 
$P(C_{n,2};s,t)$ resp. $P(T_{n,2};s,t)$ 
one obtains a valid formula for
$P(C_{n,2};t)$ resp. $P(T_{n,2};t)$. 

It is much simpler to compute the functions $P(C_{n,2};t)$ and
$P(T_{n,2};t)$ 
in the simply graded case than the corresponding functions
$P(C_{n,2};s,t)$ and $P(C_{n,2};s,t)$ in the bigraded case.
Therefore we have first computed the former and then the latter and
used the equations (\ref{veze}) as a check of correctness.

We can write the above rational functions in lowest terms as
\begin{equation} \label{poen_Ct}
P(C_{n,2};t) = \frac{ N(C_{n,2};t)}{ D(C_{n,2};t)}
\end{equation}
and
\begin{equation} \label{poen_Tt}
P(T_{n,2};t) = \frac{ N(T_{n,2};t)}{ D(T_{n,2};t)}. 
\end{equation}
{}For $n\le6$, the above denominators are products
of binomials $1-t^k$ with $k\ge1$.

Note that $N(C_{n,2};t)$ and $N(C_{n,2};t,t)$ may be different
because the numerator
$N(C_{n,2};t,t)$ and the denominator $D(C_{n,2};t,t)$ may have a
common factor. Indeed, for $n=4,5,6$ the $\gcd$ of these two
polynomials is:
\begin{eqnarray*}
&&  n=4:\quad 1+t+t^2, \\
&&  n=5:\quad (1-t)^2(1+t^2)^2, \\
&&  n=6:\quad (1+t)^2(1+t+t^2)(1-t^3)(1-t^4)^2(1-t^5)^3.
\end{eqnarray*}
Similarly, $N(T_{n,2};t,t)$ and $D(T_{n,2};t,t)$ may have a
common factor. It turns out that their gcd's are the same as
the ones listed above.

In Tables 3 and 4 we list the lowest degree numerators and
denominators for the simply graded Poincar\'{e} series
$P(C_{n,2};t)$ and $P(T_{n,2};t)$ for $n\le6$.
For the actual power series expansions of the rational functions
$P(C_{n,2};t)$ and  $P(T_{n,2};t)$ see Appendices C and D
(Tables 9 and 11), respectively.

\begin{center}

{\bf Table 3: Numerators $N(C_{n,2};t)$
and denominators $D(C_{n,2};t)$}

\begin{eqnarray*}
 N(C_{1,2};t) &=& N(C_{2,2};t)=1, \\
 D(C_{1,2};t) &=&  \left( 1-t \right) ^{2}, \\
 D(C_{2,2};t) &=&  \left( 1-t \right) ^{2} \left( 1-{t}^{2} \right) ^{3}, \\
 N(C_{3,2};t) &=& 1-{t}^{2}+{t}^{4}, \\
 D(C_{3,2};t) &=& \left( 1-t \right) ^{2} \left( 1-{t}^{2} \right) ^{4} 
\left( 1-{t}^{3} \right) ^{4}, \\
 N(C_{4,2};t) &=&  \left( 1-{t}^{2}+{t}^{4} \right)\left( 1-t-{t}^{3}+{t}^{4}
+2\,{t}^{5}+{t}^{6}-{t}^{7}-{t}^{9}+{t}^{10} \right), \\ 
 D(C_{4,2};t) &=&  \left( 1-t \right) ^{3} \left( 1-{t}^{2} \right) ^{4}
 \left( 1-{t}^{3} \right) ^{5} \left( 1-{t}^{4} \right) ^{5}, \\
 N(C_{5,2};t) &=& 1+2\,t-6\,{t}^{3}-9\,{t}^{4}+2\,{t}^{5}
+25\,{t}^{6}+38\,{t}^{7}+17\,{t}^{8}-34\,{t}^{9}-68\,{t}^{10} \\
&& -34\,{t}^{11}+73\,{t}^{12}+176\,{t}^{13}
+171\,{t}^{14}+34\,{t}^{15}-127\,{t}^{16}-156\,{t}^{17} \\
&& -2\,{t}^{18}+218\,{t}^{19}+322\,{t}^{20}+218\,{t}^{21}-\,\cdots\, 
+2\,{t}^{39}+{t}^{40}, \\
 D(C_{5,2};t) &=&  \left( 1-{t}^{2} \right) ^{6} \left( 1-{t} ^{3} \right) ^{8}
\left( 1-{t}^{4} \right) ^{6} \left( 1-{t}^{5} \right) ^{6}, \\
 N(C_{6,2};t) &=& 1-3\,t+3\,{t}^{2}-3\,{t}^{3}+3\,{t}^{4}
+4\,{t}^{5}-2\,{t}^{6}-8\,{t}^{8}-8\,{t}^{9}+11\,{t}^{10}+{t}^{11} \\
&& +56\,{t}^{12}-24\,{t}^{13}+48\,{t}^{14}-69\,{t}^{15}
-9\,{t}^{16}+2\,{t}^{17}+78\,{t}^{18}+118\,{t}^{19} \\
&& +223\,{t}^{20}+23\,{t}^{21}+158\,{t}^{22}-182\,{t}^{23}
+221\,{t}^{24}-42\,{t}^{25}+600\,{t}^{26} \\
&& +365\,{t}^{27}+633\,{t}^{28}+324\,{t}^{29}+303\,{t}^{30}
-31\,{t}^{31}+484\,{t}^{32}+178\,{t}^{33} \\
&& +1055\,{t}^{34}+ 518\,{t}^{35}
+1055\,{t}^{36}+\,\cdots\,-3\,{t}^{69}+{t}^{70}, \\
 D(C_{6,2};t) &=& 
 \left( 1-t \right) ^{5} \left( 1-{t}^{2} \right) ^{3} \left( 1-{t}^{3
} \right) ^{6} \left( 1-{t}^{4} \right) ^{9} \left( 1-{t}^{5} \right) 
^{7} \left( 1-{t}^{6} \right) ^{7}.
\end{eqnarray*}
\end{center}

\newpage
\begin{center}

{\bf Table 4: Numerators $N(T_{n,2};t)$
and denominators $D(T_{n,2};t)$}

\begin{eqnarray*}
 N(T_{1,2};t) &=& N(T_{2,2};t)=N(T_{3,2};t)=1, \\
 D(T_{1,2};t) &=&  \left( 1-t \right) ^{2}, \\
 D(T_{2,2};t) &=& 
 \left( 1-t \right) ^{4} \left( 1-{t}^{2} \right), \\
 D(T_{3,2};t) &=& 
\left( 1-t \right) ^{4} \left( 1-{t}^{2} \right) ^{4} 
\left( 1-{t}^{3} \right) ^{2}, \\
 N(T_{4,2};t) &=&  1-t+{t}^{3}+{t}^{5}-{t}^{7}+{t}^{8}, \\ 
 D(T_{4,2};t) &=& 
 \left( 1-t \right) ^{5} \left( 1-{t}^{2} \right) ^{4}
 \left( 1-{t}^{3} \right) ^{5} \left( 1-{t}^{4} \right) ^{3}, \\
 N(T_{5,2};t) &=&  1+2\,t+{t}^{2}-2\,{t}^{3}-{t}^{4}+8\,{t}^{5}+20\,{t}^{6}
+24\,{t}^{7}+18\,{t}^{8}+12\,{t}^{9} \\
&& +20\,{t}^{10}+44\,{t}^{11}+76\,{t}^{12}+94\,{t}^{13}
+85\,{t}^{14}+58\,{t}^{15}+44\,{t}^{16} \\
&& +58\,{t}^{17}+ \,\cdots\, +2\,{t}^{31}+{t}^{32}, \\
 D(T_{5,2};t) &=&  \left( 1-t \right) ^{2}
 \left( 1-{t}^{2} \right) ^{6} \left( 1-{t} ^{3} \right) ^{8}
\left( 1-{t}^{4} \right) ^{6} \left( 1-{t}^{5} \right) ^{4}, \\
 N(T_{6,2};t) &=&  1-3\,t+4\,{t}^{2}-4\,{t}^{3}+4\,{t}^{4}+3\,{t}^{5}
-6\,{t}^{6}+11\,{t}^{7}-12\,{t}^{8}+12\,{t}^{9} \\
&& +12\,{t}^{10}+{t}^{11}+55\,{t}^{12}-22\,{t}^{13}
+82\,{t}^{14}+77\,{t}^{16}+119\,{t}^{17}+84\,{t}^{18} \\
&& +234\,{t}^{19}+160\,{t}^{20}+227\,{t}^{21}+312\,{t}^{22}
+207\,{t}^{23}+507\,{t}^{24}+301\,{t}^{25} \\
&& +612\,{t}^{26}+469\,{t}^{27}+517\,{t}^{28}+593\,{t}^{29}
+426\,{t}^{30}+593\,{t}^{31} +\cdots \\
&& -3\,{t}^{59}+{t}^{60}, \\
 D(T_{6,2};t) &=& 
\left( 1-t \right) ^{7} \left( 1-{t}^{2} \right) ^{3} \left( 1-{t}^{3}
 \right) ^{6} \left( 1-{t}^{4} \right) ^{9} \left( 1-{t}^{5}
 \right) ^{7} \left( 1-{t}^{6} \right) ^{5}. \\
\end{eqnarray*}
\end{center}

It is often desirable to rewrite the formulae
(\ref{poen_Ct}) and (\ref{poen_Tt}) as
\begin{equation} \label{poen_Ct*}
P(C_{n,2};t) = \frac{ N^*(C_{n,2};t)}{ D^*(C_{n,2};t)}
\end{equation}
and
\begin{equation} \label{poen_Tt*}
P(T_{n,2};t) = \frac{ N^*(T_{n,2};t)}{ D^*(T_{n,2};t)}, 
\end{equation}
where the numerators have nonnegative integral coefficients,
and the denominators are still products of binomials $1-t^k$
with $k\ge1$. Such forms exist but are not unique.
The most interesting ones arise from a choice of a homogeneous
system of parameters (HSOP) for the algebra $C_{n,2}$ in which case
the multiset of the exponents $k$ is the same as the one made up
from degrees of the polynomials in the HSOP.
We give several such formulae in
Tables 5 and 6 when they differ from those in 
 (\ref{poen_Ct}) and (\ref{poen_Tt}).
Those in Table 5 do arise from an HSOP except possibly the one
for $C_{5,2}$.

\newpage
\begin{center}

{\bf Table 5: Numerators $N^*(C_{n,2};t)$
and denominators $D^*(C_{n,2};t)$}

\begin{eqnarray*}
 N^*(C_{3,2};t) &=&  \left( 1+{t}^{2} \right) N(C_{3,2};t) = 1+{t}^{6}, \\
 D^*(C_{3,2};t) &=&  \left( 1-t \right) ^{2} \left( 1-{t}^{2} \right) ^{3} 
\left( 1-{t}^{3} \right) ^{4} \left( 1-{t}^{4} \right), \\
 N^*(C_{4,2};t) &=&  \left( 1+{t}^{2} \right) \left( 1+t+{t}^{2} \right) 
\left( 1+{t}^{3} \right) ^{2} N(C_{4,2};t) \\
&=&  \left( 1+{t}^{6} \right)
\left( 1+2\,{t}^{5}+{t}^{6}+2\,{t}^{7}+4\,{t}^{8}+4\,{t}^{9}
+4\,{t}^{10}+\cdots+{t}^{18} \right), \\
&=&  1+2\,{t}^{5}+2\,{t}^{6}+2\,{t}^{7}+4\,{t}^{8}+4\,{t}^{9}
+4\,{t}^{10}+4\,{t}^{11}+2\,{t}^{12} \\
&& +4\,{t}^{13}+4\,{t}^{14}+4\,{t}^{15}+4\,{t}^{16}
+2\,{t}^{17}+2\,{t}^{18}+2\,{t}^{19}+{t}^{24}, \\
 D^*(C_{4,2};t) &=& 
 \left( 1-t \right) ^{2} \left( 1-{t}^{2} \right) ^{3}
 \left( 1-{t}^{3} \right) ^{4} \left( 1-{t}^{4} \right) ^{6}
 \left( 1-{t}^{6} \right) ^{2}, \\ 
 N^*(C_{5,2};t) &=&  \left( 1-t+{t}^{2} \right) ^{2} 
\left( 1+{t}^{3} \right) ^{2} \left( 1+{t}^{2}+{t}^{4} \right) 
N(C_{5,2};t) \\ 
&=&  1+2\,{t}^{5}+2\,{t}^{6}+8\,{t}^{7}+13\,{t}^{8}
+16\,{t}^{9}+25\,{t}^{10}+28\,{t}^{11}+46\,{t}^{12} \\
&& +58\,{t}^{13}+85\,{t}^{14}+132\,{t}^{15}+
172\,{t}^{16}+232\,{t}^{17}+282\,{t}^{18}+346\,{t}^{19} \\
&& +404\,{t}^{20}+444\,{t}^{21}
+518\,{t}^{22}+570\,{t}^{23}+633\,{t}^{24}+684\,{t}^{25} \\
&& +711\,{t}^{26}+744\,{t}^{27}+711\,{t}^{28}
+684\,{t}^{29}+ \cdots +2\,{t}^{49}+{t}^{54}, \\
 D^*(C_{5,2};t) &=&  \left( 1-t \right) ^{2}
 \left( 1-{t}^{2} \right) ^{3} \left( 1-{t} ^{3} \right) ^{4}
\left( 1-{t}^{4} \right) ^{6} \left( 1-{t}^{5} \right) ^{6}
\left( 1-{t}^{6} \right) ^{5}. \\
\end{eqnarray*}
\end{center}

\begin{center}

{\bf Table 6: Numerators $N^*(T_{n,2};t)$
and denominators $D^*(T_{n,2};t)$}

\begin{eqnarray*}
 N^*(T_{4,2};t) &=&  \left( 1+t+{t}^{2} \right) N(T_{4,2};t)
= 1+{t}^{4}+2\,{t}^{5}+{t}^{6}+{t}^{10}, \\ 
 D^*(T_{4,2};t) &=&  \left( 1-t \right) ^{4} \left( 1-{t}^{2} \right) ^{4}
 \left( 1-{t}^{3} \right) ^{6} \left( 1-{t}^{4} \right) ^{3}, \\
 N^*(T_{5,2};t) &=&  \left( 1+t \right) \left( 1+{t}^{2} \right)
 N(T_{5,2};t) \\
 &=&  1+3\,t+4\,{t}^{2}+2\,{t}^{3}+6\,{t}^{5}+25\,{t}^{6}
+51\,{t}^{7}+70\,{t}^{8}+74\,{t}^{9} \\
&& +74\,{t}^{10}+94\,{t}^{11}+152\,{t}^{12}+234\,{t}^{13}
+299\,{t}^{14}+313\,{t}^{15} \\
&& +281\,{t}^{16}+245\,{t}^{17}+245\,{t}^{18}
+281\,{t}^{19}+\,\cdots\, +3\,{t}^{34}+{t}^{35}, \\ 
 D^*(T_{5,2};t) &=&  \left( 1-t \right) \left( 1-{t}^{2} \right) ^{6}
 \left( 1-{t}^{3} \right) ^{8} \left( 1-{t}^{4} \right) ^{7}
 \left( 1-{t}^{5} \right) ^{4}.
\end{eqnarray*}
\end{center}

\section{ Generators of the algebra $C_{4,2}$ }
\label{gen-c42}

Teranishi \cite{T1} has constructed an HSOP for the pure trace
algebra $C_{4,2}$. Here is his result:

\begin{theorem} \label{Ter-1}
Let $x$ and $y$ be generic $4\times4$ matrices. 
Then the traces of the 17 matrices (which we arrange according
to their degrees) 
\begin{eqnarray*}
&& x,\, y, \\
&& x^2,\, xy,\, y^2, \\
&& x^3,\, x^2y,\, xy^2,\, y^3, \\
&& x^4,\, x^3y,\, x^2y^2,\, xy^3,\, y^4,\, xyxy, \\
&& (x^2y)^2,\, (y^2x)^2
\end{eqnarray*}
form an HSOP of the algebra $C_{4,2}$.
\end{theorem}

However, he did not compute a minimal set of generators for $C_{4,2}$.

The group $\GL_2$ acts (via its standard representation)
on the 2-dimensional space spanned by the
matrices $x$ and $y$. This action induces an action on the
algebra $C_{4,2}$, which was investigated by Drensky and Sadikova
in their recent paper \cite{DS}.
They show that there is a minimal set of generators whose span
is a semisimple graded $\GL_2$-submodule of $C_{4,2}$.
Moreover they determine the structure of this module.
Its Poincar\'{e} polynomial is
\[
2t+3t^2+4t^3+6t^4+2t^5+4t^6+2t^7+4t^8+4t^9+t^{10}.
\]
Hence, a minimal generating set consists of 32 polynomials.
In their paper they do not list explicitly such a generating set.

We have computed a minimal set of generators of $C_{4,2}$
independently of the above mentioned paper:

\begin{theorem} \label{Gen-1}
The 17 traces mentioned in Theorem \ref{Ter-1} together with the
traces of the following 15 matrices (arranged by their degrees)
\begin{eqnarray*}
&& x^3y^2,\, y^3x^2, \\
&& x^2y^2xy,\, y^2x^2yx, \\
&& x^3y^2xy,\, y^3x^2yx, \\
&& x^3y^2x^2y,\, y^3x^2y^2x,\, x^3y^3xy,\, y^3x^3yx, \\
&& x^3yx^2yxy,\, x^2y^2xyx^2y,\, y^2x^2yxy^2x,\, y^3xy^2xyx, \\
&& x^3y^3x^2y^2
\end{eqnarray*}
form a minimal set of generators of $C_{4,2}$.
\end{theorem}

\begin{proof}
Our proof is computational. We start with the system of parameters
from Teranishi's theorem and consider it as the
first approximation to the genuine generating set which we
want to construct. Then we compare the Poincar\'{e} series
of $C_{4,2}$ with that of its subalgebra
generated by this HSOP. The difference of the former and the latter
is a series with nonnegative integer coefficients. We
find the first nonzero term, say $ct^d$. This means that we have
to enlarge our incomplete generating set by additional $c$
generators of degree $d$. We then select $c$ words in $x$ and $y$
of length $d$, whose traces provide these additional generators.

We repeat this procedure with the enlarged set of generators,
and continue repeating it until we reach the space of invariants
of degree 10. After adding the single additional generator in
degree 10, we are certain that we have indeed found the full
set of generators. This is a consequence of the well known fact
\cite[Part A, Theorem 6.1.6]{DF}, originally due to
Procesi \cite{CPr},
which establishes the connection between the Nagata--Higman
theorem and the invariant theory of generic matrices.
\end{proof}

\begin{remark} \label{primedba}
The words $y^3x^2yx$, $y^3x^2y^2x$, $y^2x^2yxy^2x$, $y^3xy^2xyx$ in
Theorem \ref{Gen-1} can be replaced with $xyx^2y^3$, $xy^2x^2y^3$,
$xy^2xyx^2y^2$, $xyxy^2xy^3$, respectively.
\end{remark}

Note that the total number of generators listed in the theorem
is 32, and that this number as well as the degrees of
these generators are in agreement with the result of Drensky
and Sadikova mentioned above.

Let us say that a set $W$ of words in two noncommuting indeterminates
$x,y$ is a {\em test set} for $M_n(\bC)$ if it has the following
property: Two matrices $a,b\in M_n(\bC)$ are unitarily similar iff
$\tr(w(a,a^*))=\tr(w(b,b^*))$ holds for all $w\in W$. Such a test
set $W$ is {\em minimal} if no proper subset of $W$ is a test set.

As a consequence of the above theorem, we can state the following
criterion for unitary similarity of two $4\times4$ complex
matrices.

\begin{theorem} \label{Un-sim}
The following 20 words form a test set for $M_4(\bC)$:
\begin{eqnarray*}
&& x;\, x^2,\, xy;\, x^3,\, x^2y;\,
x^4,\, x^3y,\, x^2y^2,\, xyxy;\, x^3y^2;\, (x^2y)^2,\, x^2y^2xy,\, y^2x^2yx; \\
&& x^3y^2xy;\, x^3y^2x^2y,\, x^3y^3xy,\, y^3x^3yx;\,
x^3yx^2yxy,\, x^2y^2xyx^2y;\, x^3y^3x^2y^2.
\end{eqnarray*}

\end{theorem}

\begin{proof}
First apply the above remark. Then
among the 32 words $w(x,y)$ whose traces generate the algebra
$C_{4,2}$, given in Theorem \ref{Gen-1}, there are 12 pairs
$\{w_1,w_2\}$ such that, for any $4\times 4$ complex matrix
$a$, $\tr(w_1(a,a^*))$ and $\tr(w_2(a,a^*))$ are complex conjugates.
For instance, $\{x,y\}$, $\{x^2,y^2\}$  and $\{x^3y^2,y^3x^2\}$
are such pairs.
(For the non-paired words, such as $w(x,y)=xyxy$, the trace of
$w(a,a^*)$ is always real.)
By dropping one of the words from each of these pairs,
we obtain the test set in the theorem.
\end{proof}

\section{ Generators of the algebra $C_{5,2}$ }
\label{gen-c52}

Let $M_n(0)$ be the subspace of $M_n$ consisting of matrices
of trace 0 and let
$C_{n,2}(0)=K[M_n(0)^2]^{\GL_n}$
be the corresponding algebra of $\GL_n$-invariant polynomial
functions on the direct product
$M_n(0)^2=M_n(0)\times M_n(0)$.
Then one has an isomorphism of $\bZ^2$-graded algebras
\begin{equation} \label{tens-pr}
 C_{n,2}\cong K[u,v]\otimes C_{n,2}(0),
\end{equation}
where $K[u,v]$ is the polynomial algebra in two variables $u$ and $v$
(see e.g. \cite[Section 5]{CPr}).
Consequently, the problem of constructing a minimal set of generators of
$C_{n,2}$ reduces to the same problem for $C_{n,2}(0)$.
In the remainder of this section we shall assume that
$x,y\in M_n(0)$. (More precisely, the entries of $x$ and $y$ will be treated,
via restriction, as linear functions on the subspace $M_n(0)$ of $M_n$.)

In view of the above isomorphism, we have the following obvious relation
between the Poincar\'{e} series of these algebras:
\[ P(C_{n,2};t)=\frac{ P(C_{n,2}(0);t) }{ (1-t)^2 }. \]

\begin{theorem} \label{gen-C52}
The algebra $C_{5,2}(0)$ has a minimal generating set $\pP$
consisting of 171 bihomogeneous polynomials.
$\pP$ is the disjoint union of four subsets:
$\pP_s$, $\pP_k$, $\pP_d$ and $\pP'_d$ with cardinals 5, 4, 81
and 81, respectively. The polynomials $f(x,y)\in\pP_s$ are symmetric,
(i.e., satisfy $f(y,x)=f(x,y)$), while those in $\pP_k$ are skew-symmetric.
The polynomials in $\pP_d$ are neither symmetric nor skew-symmetric.
There is a bijection $\pP_d\to\pP'_d$ given by $f(x,y)\to f(y,x)$.
The sets $\pP_s$, $\pP_k$ and $\pP_d$ are given below.
\end{theorem}

$\pP_s$ consists of the traces of the following 5 words
in $x$ and $y$:
\[ xy,\, x^2y^2,\, xyxy,\, x^3y^3,\, x^4y^4. \]
$\pP_k$ consists of the traces of the following 4 matrices:
\begin{eqnarray*}
&&  x^3y^3x^2y^2-y^3x^3y^2x^2,\, x^2yxy^2xyxy-y^2xyx^2yxyx, \\
&& x^3y^2xyxy^2xy-y^3x^2yxyx^2yx,\, x^4y^4x^3y^3-y^4x^4y^3x^3.
\end{eqnarray*}
Finally, the set $\pP_d$ consists of the traces of the following
81 words:
\begin{eqnarray*}
&& x^2; \\
&& x^3,\, x^2y; \\
&& x^4,\, x^3y; \\
&& x^5,\, x^4y,\, x^3y^2,\, x^2yxy; \\
&& x^4y^2,\, (x^2y)^2,\, x^2y^2xy; \\
&& x^4y^3,\, x^4yxy,\, x^3y^2xy,\, x^3yxy^2; \\
&& x^4y^2xy,\, x^4yx^2y,\, x^4yxy^2,\, x^3y^3xy,\, x^3y^2x^2y,\, x^2y^2xyxy; \\
&& x^4y^2x^2y,\, x^4y^2xy^2,\, x^4yx^2y^2,\, x^3y^3x^2y,\, x^3y^2xyxy,\,
x^3yxyx^2y,\, x^3yxyxy^2, \\
&& \quad x^2y^2xyx^2y; \\
&& x^4y^4xy,\, x^4y^3xy^2,\, x^4y^2x^3y,\, x^4y^2x^2y^2,\, x^4y^2xyxy,\,
x^4yx^2yxy,\, x^4yxy^2xy, \\
&& \quad x^3y^3xyxy,\, x^3y^2x^2yxy,\, x^3yx^2yxy^2; \\
&& x^4y^4x^2y,\, x^4y^4xy^2,\, x^4y^3x^3y,\, x^4y^3x^2y^2,\, x^4y^3xyxy,\, 
x^4y^2x^2yxy,\, x^4y^2xyxy^2, \\
&& \quad x^4yx^3yxy,\, x^4yx^2y^2xy,\, x^4yx^2yxy^2,\, x^3y^2x^2y^2xy,\, 
x^3yxyxy^2xy; \\
&& x^4y^4x^3y,\, x^4y^4x^2y^2,\, x^4y^3x^3y^2,\, x^4y^3x^2yxy,\,
x^4y^3xy^2xy,\, x^4y^2x^3yxy, \\
&& \quad x^4y^2x^2yx^2y,\, x^4y^2xy^2x^2y,\, x^4y^2xyx^2y^2,\,
x^4yx^3yx^2y,\, x^4yx^3yxy^2, \\
&& \quad x^4yxyxy^2xy,\, x^3y^3x^2yxy^2; \\
&& x^4y^4x^3y^2,\, x^4y^4x^2yxy,\, x^4y^4xy^2xy,\, x^4y^3x^3yxy,\,
x^4y^3x^2y^2xy,\, x^4y^3(x^2y)^2, \\
&& \quad x^4y^2x^3yx^2y,\, x^4(y^2x^2)^2y,\, x^4y^2xyxy^2xy; \\
&& x^4y^4x^3yxy,\, x^4y^4x^2y^2xy,\, x^4y^3x^3y^2xy,\, x^4y^3x^2yx^3y; \\
&& x^4y^4x^3y^2xy,\, x^4y^3x^3y^2x^2y,\, x^4y^2x^3yx^2yxy.
\end{eqnarray*}

\begin{proof}
Shestakov and Zhukavets \cite{SZ} have shown that any
2-generated associative algebra
(non-unital and over a field of characteristic 0)
which satisfies the identity $x^5=0$, also satisfies the
identity $x_1x_2\cdots x_{15}=0$.
Since we are working with only two generic matrices, by invoking
a theorem of Procesi \cite[Theorem 3.2]{CPr}, we conclude that the
algebra $C_{5,2}$ (and $C_{5,2}(0)$) is generated by polynomials of
degree at most 15.

Let $\pA$ denote the unital subalgebra of $C_{5,2}(0)$ generated
by 171 polynomials in our set $\pP$. We have verified, using a
computer, that for each degree $d\le 15$ the homogeneous component,
$\pA_d$, of $\pA$ of degree $d$ has the dimension equal to the coefficient
of $t^d$ in the Poincar\'{e} series of the algebra $C_{5,2}(0)$.
For the reader's convenience, these coefficients are given in Appendix C.
As an additional check, we have computed the dimension of $\pA_{16}$
and verified that it is indeed equal to $17338$.
\end{proof}

\begin{remark} 
In view of (\ref{tens-pr}), it follows that a minimal generating set
of the algebra $C_{5,2}$ has cardinal 173.
\end{remark}

\begin{remark} 
One can modify the minimal generating set $\pP$ by replacing each
of the four generators $\tr(w(x,y)-w(y,x))\in\pP_k$ with
$\tr(w(x,y))$. The only reason for our choice was to make $\pP$
stable (up to sign) under the involution
which interchanges $x$ and $y$.
\end{remark}

\begin{remark}
By using this theorem, one can now construct easily a test set
for $M_5(\bC)$ of cardinal $5+4+82=91$.
\end{remark}

\section{ Poincar\'{e} series for invariants of $\GL_n\times\SL_d$ }
\label{veca-gr}

The vector space $M_n^d$ can be identified with the tensor product
$M_n\otimes K^d$ of $M_n$ and the $d$-dimensional space $K^d$. The action
of $\GL_n$ on $M_n^d$ corresponds to its action on this tensor product
given by $a\cdot(x\otimes v)=axa^{-1}\otimes v$. We can now view
$M_n\otimes K^d$ as a module for the direct product $\GL_n\times\GL_d$
by letting $\GL_d$ act on $K^d$ by multiplication.
Denote by $C_{n,d}^\#$ the subalgebra  of $K[M_n^d]$ consisting of
$\GL_n\times\SL_d$-invariant functions. This is a subalgebra of $C_{n,d}$.

Our objective in this section is to record the Poincar\'{e} series
of the algebras $C_{n,2}^\#$ for $n\le5$.
We shall write these Poincar\'{e} series as rational functions
in lowest terms 
\[
P(C_{n,2}^\#;t)=\frac{N(C_{n,2}^\#;t)}{D(C_{n,2}^\#;t)},
\]
with numerator and denominator normalized so to take value 1 at $t=0$.
These numerators are again palindromic.

\begin{center}

{\bf Table 7: Numerators $N(C_{n,2}^\#;t)$
and denominators $D(C_{n,2}^\#;t)$}

\begin{eqnarray*}
 N(C_{2,2}^\#;t) &=& 1, \\
 D(C_{2,2}^\#;t) &=&  \left( 1-{t}^{4} \right) ^{2}, \\
 N(C_{3,2}^\#;t) &=& 1+3\,{t}^{8}+{t}^{10}+3\,{t}^{12}+{t}^{20}, \\
 D(C_{3,2}^\#;t) &=&  \left( 1-{t}^{4} \right) ^{3}
 \left( 1-{t}^{6} \right) ^{3} \left( 1-{t}^{8} \right), \\
 N(C_{4,2}^\#;t) &=& 1+10\,{t}^{8}+12\,{t}^{10}+38\,{t}^{12}
+46\,{t}^{14}+93\,{t}^{16}+131\,{t}^{18}+235\,{t}^{20} \\
&& +299\,{t}^{22}+473\,{t}^{24}+560\,{t}^{26}
+714\,{t}^{28}+761\,{t}^{30}+876\,{t}^{32} \\
&& +830\,{t}^{34}+876\,{t}^{36}+\cdots +12\,{t}^{58}+10\,{t}^{60}+{t}^{68}, \\
 D(C_{4,2}^\#;t) &=& \left( 1-{t}^{4} \right) ^{3} \left( 1-{t}^{6} \right) ^{4}
 \left( 1-{t}^{8} \right) ^{4} \left( 1-{t}^{10} \right) ^{2}
\left( 1-{t}^{12} \right), \\
 N(C_{5,2}^\#;t) &=& 1+{t}^{2}+{t}^{4}+{t}^{6}+14\,{t}^{8}
+41\,{t}^{10}+135\,{t}^{12}+329\,{t}^{14}+842\,{t}^{16} \\
&& +1980\,{t}^{18}+4677\,{t}^{20}
+10386\,{t}^{22}+22654\,{t}^{24}+47093\,{t}^{26} \\
&& +94970\,{t}^{28}+184182\,{t}^{30}+346523\,{t}^{32}
+629769\,{t}^{34}+1111589\,{t}^{36} \\
&& +1902191\,{t}^{38}+3165521\,{t}^{40}+5120359\,{t}^{42}
+8066607\,{t}^{44} \\
&& +12376177\,{t}^{46}+18520117\,{t}^{48}+27035364\,{t}^{50}
+38541637\,{t}^{52} \\
&& +53673328\,{t}^{54}+73078953\,{t}^{56}
+97307914\,{t}^{58}+126802726\,{t}^{60} \\
&& +161749890\,{t}^{62}+202084191\,{t}^{64}+247338162\,{t}^{66}
+296695937\,{t}^{68} \\
&& +348874713\,{t}^{70}+402270954\,{t}^{72}+454898759\,{t}^{74}
+504632564\,{t}^{76} \\
&& +549206297\,{t}^{78}+586521387\,{t}^{80}+614654835\,{t}^{82}
+632178319\,{t}^{84} \\
&& +638112785\,{t}^{86} +632178319\,{t}^{88}+ \cdots +14\,{t}^{164}
+{t}^{166}+{t}^{168} \\
&& +{t}^{170}+{t}^{172}, \\
 D(C_{5,2}^\#;t) &=&   \left( 1+{t}^{2}+{t}^{4} \right)
 \left( 1-{t}^{4} \right) ^{3} \left( 1-{t}^{6} \right) ^{3}
\left( 1-{t}^{8} \right) ^{5} \left( 1-{t}^{10} \right) ^{5}
 \left( 1-{t}^{12} \right) ^{3} \cdot \\
&& \left( 1-{t}^{14} \right) ^{2}
\left( 1-{t}^{16} \right)  \left( 1-{t}^{18} \right).
\end{eqnarray*}

\end{center}

These rational functions were again computed by using the
Molien--Weyl formula. In this case the formula is more complicated:
\[ P(C_{n,2}^\#;t) = \frac{1}{2\pi i} \int_{|y|=1}
\frac{(1-y^2)\psi(y)}{(1-ty)^n(1-\frac{t}{y})^n} \frac{{\rm d}y}{y}, \]
where
\[ \psi(y)=\frac{1}{(2\pi i)^{n-1}} \int_{|x_1|=1}\cdots\int_{|x_{n-1}|=1}
\prod_{1\le k\le r\le n-1} \frac{1-x_kx_{k+1}\cdots x_r}{\psi_{k,r}}
\frac{{\rm d}x_{n-1}}{x_{n-1}}\cdots \frac{{\rm d}x_1}{x_1}, \]
\begin{eqnarray*}
\psi_{k,r} &=& (1-tyx_k x_{k+1}\cdots x_r)(1-ty^{-1}x_k x_{k+1}\cdots x_r)\cdot \\
&& (1-ty(x_k x_{k+1}\cdots x_r)^{-1})(1-t(yx_k x_{k+1}\cdots x_r)^{-1})
\end{eqnarray*}
and the integration is performed over the unit circles assuming that $|t|<1$.

The computations were easy for $n\le4$ but very hard (lasting several
days) for $n=5$. We close this section by listing the terms
of degree $<26$ in the Taylor series of these rational functions.

\begin{eqnarray*}
 P(C_{2,2}^\#;t) &=& 
1+2\,{t}^{4}+3\,{t}^{8}+4\,{t}^{12}+5\,{t}^{16}+6\,{t}^{20}+7\,{t}^{24}
 +\cdots, \\
 P(C_{3,2}^\#;t) &=& 1+3\,{t}^{4}+3\,{t}^{6}+10\,{t}^{8}+10\,{t}^{10}
+31\,{t}^{12}+33\,{t}^{14}+73\,{t}^{16}+92\,{t}^{18} \\
&& +164\,{t}^{20}+205\,{t}^{22}+344\,{t}^{24} +\cdots, \\
 P(C_{4,2}^\#;t) &=& 1+3\,{t}^{4}+4\,{t}^{6}+20\,{t}^{8}+26\,{t}^{10}
+101\,{t}^{12}+168\,{t}^{14}+445\,{t}^{16} \\
&& +813\,{t}^{18}+1804\,{t}^{20}+3246\,{t}^{22}+6527\,{t}^{24} +\cdots, \\
 P(C_{5,2}^\#;t) &=&   1+3\,{t}^{4}+4\,{t}^{6}+24\,{t}^{8}+44\,{t}^{10}
+171\,{t}^{12}+388\,{t}^{14}+1166\,{t}^{16} \\
&& +2808\,{t}^{18}+7344\,{t}^{20}+17281\,{t}^{22}+41569\,{t}^{24} +\cdots.
\end{eqnarray*}

Since the algebra $C_{5,2}^\#$ is the algebra of $\SL_2$-invariants
in $C_{5,2}$, the coefficients of $t^{2k}$ in the above Taylor
series must be the same as the coefficients of the Schur functions $f_{k,k}$
in the formula for $P(C_{5,2};t)$ displayed in Appendix B.
It is easy to check that this is indeed the case for $k\le12$, which
gives a further confirmation of our formula for $P(C_{5,2}^\#;t)$.

\section{ Poincar\'{e} series for invariants of $\GL_n\times\Delta_{d-1}$ }
\label{manja-gr}

Let us restrict the action of $\GL_n\times\SL_d$ on
$M_n^d$ to its subgroup $\GL_n\times\Delta_{d-1}$, where 
$\Delta_{d-1}$ is the maximal torus of $\SL_d$
consisting of diagonal matrices.
Denote by $C_{n,d}^\bullet$ the subalgebra of $K[M_n^d]$ consisting
of $\GL_n\times\Delta_{n-1}$-invariant polynomial functions.
This is a subalgebra of $C_{n,d}$.

In this section we record the Poincar\'{e} series
of $C_{n,2}^\bullet$ for $n\le5$.
We shall write these rational functions in lowest terms as
\[
P(C_{n,2}^\bullet;t)=\frac{N(C_{n,2}^\bullet;t)}{D(C_{n,2}^\bullet;t)},
\]
with numerator and denominator normalized so to take value 1 at $t=0$.
These numerators are again palindromic.

\begin{center}

{\bf Table 8: Numerators $N(C_{n,2}^\bullet;t)$
and denominators $D(C_{n,2}^\bullet;t)$}

\begin{eqnarray*}
&& N(C_{1,2}^\bullet;t)=1, \\
&& D(C_{1,2}^\bullet;t)=1-{t}^{2}, \\
&& N(C_{2,2}^\bullet;t)=1+{t}^{4}, \\
&& D(C_{2,2}^\bullet;t)= \left( 1-{t}^{2} \right) ^{2}
\left( 1-{t}^{4} \right) ^{2}, \\
&& N(C_{3,2}^\bullet;t)=1+3\,{t}^{4}+6\,{t}^{6}+9\,{t}^{8}+6\,{t}^{10}
+12\,{t}^{12}+6\,{t}^{14}+ \cdots +{t}^{24}, \\
&& D(C_{3,2}^\bullet;t)= 
 \left( 1-{t}^{2} \right) ^{2} \left( 1-{t}^{4} \right) ^{3} \left( 1-
{t}^{6} \right) ^{3} \left( 1-{t}^{8} \right), \\
&& N(C_{4,2}^\bullet;t)= \\
&& 1+4\,{t}^{4}+12\,{t}^{6}+36\,{t}^{8}+68\,{t}^{10}
+171\,{t}^{12}+316\,{t}^{14}+639\,{t}^{16}+1096\,{t}^{18} \\
&& +1096\,{t}^{18}+1849\,{t}^{20}+2794\,{t}^{22}+4151\,{t}^{24}
+5546\,{t}^{26}+7229\,{t}^{28}+8700\,{t}^{30} \\
&& +10085\,{t}^{32}+10836\,{t}^{34}+11270\,{t}^{36}+10836\,{t}^{38}
+ \cdots +12\,{t}^{66}+4\,{t}^{68}+{t}^{72}, \\
&& D(C_{4,2}^\bullet;t)=  \left( 1-{t}^{2} \right) ^{2}
\left( 1-{t}^{4} \right) ^{3} 
\left( 1-{t}^{6} \right) ^{4} \left( 1-{t}^{8} \right) ^{4}
\left( 1-{t}^{10} \right) ^{2} \left( 1-{t}^{12} \right), \\
&& N(C_{5,2}^\bullet;t)= \\ 
&& 1+{t}^{2}+5\,{t}^{4}+20\,{t}^{6}+76\,{t}^{8}+227\,{t}^{10}
+692\,{t}^{12}+1933\,{t}^{14}+5307\,{t}^{16}+13752\,{t}^{18} \\
&& +34304\,{t}^{20}+81525\,{t}^{22}+186346\,{t}^{24}+408071\,{t}^{26}
+860437\,{t}^{28}+1746504\,{t}^{30} \\
&& +3421732\,{t}^{32}+6474866\,{t}^{34}+11857662\,{t}^{36}
+21033945\,{t}^{38}+36195856\,{t}^{40} \\
&& +60479854\,{t}^{42}+98242554\,{t}^{44}
+155273212\,{t}^{46}+239019423\,{t}^{48}+358621723\,{t}^{50} \\
&& +524884888\,{t}^{52}+749897456\,{t}^{54}+1046516425\,{t}^{56}
+1427383948\,{t}^{58} \\
&& +1903851664\,{t}^{60}+2484438301\,{t}^{62}
+3173436196\,{t}^{64}+3969248353\,{t}^{66} \\
&& +4863282209\,{t}^{68}+5838905156\,{t}^{70}+6871421892\,{t}^{72}
+7928353846\,{t}^{74} \\
&& +8971036674\,{t}^{76}+9956478001\,{t}^{78}+10840418189\,{t}^{80}
+11580232480\,{t}^{82} \\
&& +12138549745\,{t}^{84}+12485984964\,{t}^{86}+
12603960344\,{t}^{88}+12485984964\,{t}^{90} \\
&& + \cdots +20\,{t}^{170}+5\,{t}^{172}+{t}^{174}+{t}^{176}, \\
&& D(C_{5,2}^\bullet;t)= \\  
&& \left( 1-{t}^{2} \right)  \left( 1-{t}^{4} \right) ^{3} \left( 1-{t}^
{6} \right) ^{4} \left( 1-{t}^{8} \right) ^{5} \left( 1-{t}^{10}
 \right) ^{5} \left( 1-{t}^{12} \right) ^{3} \left( 1-{t}^{14} \right) ^{2}
\left( 1-{t}^{16} \right) \cdot \\
&& \left( 1-{t}^{18} \right).
\end{eqnarray*}

\end{center}

These rational functions were computed by using the formula:
\[ P(C_{n,2}^\bullet;t) = \frac{1}{2\pi i} \int_{|z|=1}
P(C_{n,2};tz^{-1},tz) \frac{{\rm d}z}{z}, \]
where $P(C_{n,2};s,t)$ is the bigraded Poincar\'{e} series of $C_{n,2}$
from Section \ref{z2-grad}.
We list the terms of degree $<26$ in their Taylor expansions:

\newpage
\begin{eqnarray*}
P(C_{1,2}^\bullet;t) &=& 
1+{t}^{2}+{t}^{4}+{t}^{6}+{t}^{8}+{t}^{10}+{t}^{12}+{t}^{14}+{t}^{16}
+{t}^{18}+{t}^{20}+{t}^{22} \\
&& +{t}^{24} +\cdots, \\
 P(C_{2,2}^\bullet;t) &=& 1+2\,{t}^{2}+6\,{t}^{4}+10\,{t}^{6}+19\,{t}^{8}
+28\,{t}^{10}+44\,{t}^{12}+60\,{t}^{14}+85\,{t}^{16} \\
&& +110\,{t}^{18}+146\,{t}^{20}+182\,{t}^{22}+231\,{t}^{24} +\cdots,  \\
 P(C_{3,2}^\bullet;t) &=& 
1+2\,{t}^{2}+9\,{t}^{4}+25\,{t}^{6}+66\,{t}^{8}+149\,{t}^{10}+329\,{t}
^{12}+650\,{t}^{14} \\
&& +1248\,{t}^{16}+2255\,{t}^{18}+3941\,{t}^{20}+6608\,{t}^{22}+10799\,{t}^{24} +\cdots, \\
 P(C_{4,2}^\bullet;t) &=& 
1+2\,{t}^{2}+10\,{t}^{4}+34\,{t}^{6}+116\,{t}^{8}+332\,{t}^{10}+948\,{
t}^{12}+2450\,{t}^{14} \\
&& +6126\,{t}^{16}+14426\,{t}^{18}+32746\,{t}^{20}+71100\,{t}^{22}+149402\,{t}^{24}
+\cdots,  \\
 P(C_{5,2}^\bullet;t) &=&   
1+2\,{t}^{2}+10\,{t}^{4}+37\,{t}^{6}+143\,{t}^{8}+478\,{t}^{10}+1588\,
{t}^{12}+4910\,{t}^{14} \\
&& +14748\,{t}^{16}+42235\,{t}^{18}+116910\,{t}^{20}+311478\,{t}^{22}
+803343\,{t}^{24}  \\
&& +\cdots.
\end{eqnarray*}

\section{ Conjectures }
\label{hipot}

On the basis of our computations, we propose four conjectures
about the numerators 
$N(C_{n,2};s,t)$ and $N(T_{n,2};s,t)$
and the denominators 
$D(C_{n,2};s,t)$ and $D(T_{n,2};s,t)$.
See Section \ref{z2-grad} for their definitions.

\begin{conjecture} \label{conj-1} The denominators 
$D(C_{n,2};s,t)$ and $D(T_{n,2};s,t)$ can be written as products
of binomials $1-s^a t^b$, where $a$ and $b$ are nonnegative integers.
\end{conjecture} 

\begin{conjecture} \label{conj-2} $N(C_{n,2};1,1)=N(T_{n,2};1,1)=0$
for $n\ge5$.
\end{conjecture} 

\begin{conjecture} \label{conj-3} For all $n$,
\[ (1-s)(1-t)D(C_{n,2};s,t)=(1-s^n)(1-t^n)D(T_{n,2};s,t). \]
\end{conjecture}

\begin{conjecture} \label{conj-4} For all $n$,
\[ \gcd\left( N(C_{n,2};t,t),D(C_{n,2};t,t) \right)
=\gcd\left( N(T_{n,2};t,t),D(T_{n,2};t,t) \right). \]
\end{conjecture}

All four conjectures are true for $n\le6$. The second conjecture
implies (see the proof of Proposition \ref{nema-Z2}) that
$C_{n,2}$ has no bigraded system of parameters for $n\ge5$.
As all coefficients of the numerators
$N(C^\#_{n,2};t)$ and $N(C^\bullet_{n,2};t)$
given in Sections \ref{veca-gr} and \ref{manja-gr} are
nonnegative, we propose yet another conjecture.

\begin{conjecture} \label{conj-5} For all $n$,
$N(C^\#_{n,2};t)$ and $N(C^\bullet_{n,2};t)$ have
nonnegative integer coefficients.
\end{conjecture}

The following interesting problem may have some practical applications.

\begin{problem} \label{prob-1}
Construct minimal test sets for $M_4(\bC)$ and $M_5(\bC)$.
\end{problem}

\section{ Appendix A: Bigraded Poincar\'{e} series for 
$C_{6,2}$  and $T_{6,2}$}

The numerator $N(C_{6,2};s,t)$ is a polynomial in $s$ and $t$ of
total degree 100, with leading term $(st)^{50}$,
having in total 1169 terms. By using its symmetry properties,
we can write it as follows
\[ N(C_{6,2};s,t)=f(s,t)+(st)^{50}f(s^{-1},t^{-1})-18142(st)^{25}, \]
where
\begin{eqnarray*}
f(s,t) &=& f_1(st)-2s^{17}t^{33}f_2(st^{-1})
-st^2f_3(s,t)-s^2tf_3(t,s), \\
f_1(x) &=& 1-2\,x+9\,{x}^{3}-3\,{x}^{4}-35\,{x}^{5}+27\,{x}^{6}
+137\,{x}^{7}-89\,{x}^{8}-388\,{x}^{9} \\
&& +337\,{x}^{10}+955\,{x}^{11}-704\,{x}^{12}
-2155\,{x}^{13}+1319\,{x}^{14}+4002\,{x}^{15} \\
&& -2209\,{x}^{16}-6963\,{x}^{17}+2820\,{x}^{18}
+10659\,{x}^{19}-2638\,{x}^{20}-14080\,{x}^{21} \\
&& +2275\,{x}^{22}+16918\,{x}^{23}-1114\,{x}^{24}, \\
f_2(x) &=& 10+77\,x+345\,{x}^{2}+1073\,{x}^{3}+2480\,{x}^{4}
+4519\,{x}^{5}+6700\,{x}^{6}+8497\,{x}^{7},
\end{eqnarray*}
and the polynomial $f_3$ has 276 terms:
\begin{eqnarray*}
 && f_3(s,t) = \\
 && 2+2\,t+{t}^{2}-4\,st-7\,s{t}^{2}-8\,s{t}^{3}-8\,{s}^{2}{t}^{2}
-5\,s{t}^{4}-{s}^{2}{t}^{3}-3\,s{t}^{5}+10\,{s}^{2}{t}^{4} \\
&& +21\,{s}^{3}{t}^{3}-s{t}^{6}+15\,{s}^{2}{t}^{5}+32\,{s}^{3}{t}^{4}
+15\,{s}^{2}{t}^{6}+25\,{s}^{3}{t}^{5}+25\,{s}^{4}{t}^{4}
+10\,{s}^{2}{t}^{7} \\
&& +4\,{s}^{3}{t}^{6}-18\,{s}^{4}{t}^{5}+5\,{s}^{2}{t}^{8}-13\,{s}^{3}{t}^{7}
-61\,{s}^{4}{t}^{6}-94\,{s}^{5}{t}^{5}+2\,{s}^{2}{t}^{9}
-24\,{s}^{3}{t}^{8} \\
&& -75\,{s}^{4}{t}^{7}-121\,{s}^{5}{t}^{6}-21\,{s}^{3}{t}^{9}
-61\,{s}^{4}{t}^{8}-74\,{s}^{5}{t}^{7}
-80\,{s}^{6}{t}^{6}-15\,{s}^{3}{t}^{10}-23\,{s}^{4}{t}^{9} \\
&& +10\,{s}^{5}{t}^{8}+73\,{s}^{6}{t}^{7}-6\,{s}^{3}{t}^{11}+6\,{s}
^{4}{t}^{10}+90\,{s}^{5}{t}^{9}+201\,{s}^{6}{t}^{8}
+290\,{s}^{7}{t}^{7}-2\,{s}^{3}{t}^{12} \\
&& +23\,{s}^{4}{t}^{11}+111\,{s}^{5}{t}^{10}+241\,{s}^{6}{t}^{9}
+326\,{s}^{7}{t}^{8}+21\,{s}^{4}{t}^{12}+91\,{s}^{5}{t}^{11}
+152\,{s}^{6}{t}^{10} \\
&& +173\,{s}^{7}{t}^{9}+170\,{s}^{8}{t}^{8}+13\,{s}^{4}{t}^{13}
+43\,{s}^{5}{t}^{12}+22\,{s}^{6}{t}^{11}-111\,{s}^{7}{t}^{10}
-233\,{s}^{8}{t}^{9} \\
&& +5\,{s}^{4}{t}^{14}+2\,{s}^{5}{t}^{13}-90\,{s}^{6}{t}^{12}
-318\,{s}^{7}{t}^{11}-595\,{s}^{8}{t}^{10}-808\,{s}^{9}{t}^{9}
+{s}^{4}{t}^{15} \\
&& -14\,{s}^{5}{t}^{14}-125\,{s}^{6}{t}^{13} -367\,{s}^{7}{t}^{12}
-629\,{s}^{8}{t}^{11}-902\,{s}^{9}{t}^{10}-14\,{s}^{5}{t}^{15}
-94\,{s}^{6}{t}^{14} \\
&& -248\,{s}^{7}{t}^{13}-357\,{s}^{8}{t}^{12}-430\,{s}^{9}{t}^{11}
-420\,{s}^{10}{t}^{10}-7\,{s}^{5}{t}^{16}-45\,{s}^{6}{t}^{15}
-70\,{s}^{7}{t}^{14} \\
&& +48\,{s}^{8}{t}^{13}+276\,{s}^{9}{t}^{12}+608\,{s}^{10}{t}^{11}
-2\,{s}^{5}{t}^{17}-7\,{s}^{6}{t}^{16}+56\,{s}^{7}{t}^{15}
+347\,{s}^{8}{t}^{14} \\
&& +812\,{s}^{9}{t}^{13}+1416\,{s}^{10}{t}^{12}+1807\,{s}^{11}{t}^{11}
+7\,{s}^{6}{t}^{17}+95\,{s}^{7}{t}^{16}+401\,{s}^{8}{t}^{15} \\
&& +887\,{s}^{9}{t}^{14}+1500\,{s}^{10}{t}^{13}+1947\,{s}^{11}{t}^{12}
+6\,{s}^{6}{t}^{18}+70\,{s}^{7}{t}^{17}+274\,{s}^{8}{t}^{16} \\
&& +535\,{s}^{9}{t}^{15}+830\,{s}^{10}{t}^{14}+943\,{s}^{11}{t}^{13}
+984\,{s}^{12}{t}^{12}+2\,{s}^{6}{t}^{19}+29\,{s}^{7}{t}^{18} \\
&& +101\,{s}^{8}{t}^{17}+69\,{s}^{9}{t}^{16}-131\,{s}^{10}{t}^{15}
-593\,{s}^{11}{t}^{14}-1082\,{s}^{12}{t}^{13}+5\,{s}^{7}{t}^{19} \\
&& -15\,{s}^{8}{t}^{18}-238\,{s}^{9}{t}^{17}-780\,{s}^{10}{t}^{16}
-1717\,{s}^{11}{t}^{15}-2716\,{s}^{12}{t}^{14}-3416\,{s}^{13}{t}^{13} \\
&& -2\,{s}^{7}{t}^{20}-44\,{s}^{8}{t}^{19}-299\,{s}^{9}{t}^{18}
-857\,{s}^{10}{t}^{17}-1803\,{s}^{11}{t}^{16}-2859\,{s}^{12}{t}^{15} \\
&& -3687\,{s}^{13}{t}^{14}-{s}^{7}{t}^{21}-30\,{s}^{8}{t}^{20}
-193\,{s}^{9}{t}^{19}-529\,{s}^{10}{t}^{18}-1016\,{s}^{11}{t}^{17} \\
&& -1488\,{s}^{12}{t}^{16}-1789\,{s}^{13}{t}^{15}-1873\,{s}^{14}{t}^{14}
-9\,{s}^{8}{t}^{21}-69\,{s}^{9}{t}^{20}-110\,{s}^{10}{t}^{19} \\
&& -31\,{s}^{11}{t}^{18}+396\,{s}^{12}{t}^{17}+1118\,{s}^{13}{t}^{16}
+1866\,{s}^{14}{t}^{15}-{s}^{8}{t}^{22}-{s}^{9}{t}^{21}
+122\,{s}^{10}{t}^{20} \\
&& +582\,{s}^{11}{t}^{19}+1589\,{s}^{12}{t}^{18}+3168\,{s}^{13}{t}^{17}
+4811\,{s}^{14}{t}^{16}+5888\,{s}^{15}{t}^{15}+12\,{s}^{9}{t}^{22} \\
&& +155\,{s}^{10}{t}^{21}+630\,{s}^{11}{t}^{20}+1669\,{s}^{12}{t}^{19}
+3250\,{s}^{13}{t}^{18}+5064\,{s}^{14}{t}^{17}+6457\,{s}^{15}{t}^{16} \\
&& +5\,{s}^{9}{t}^{23}+87\,{s}^{10}{t}^{22}+362\,{s}^{11}{t}^{21}
+928\,{s}^{12}{t}^{20}+1783\,{s}^{13}{t}^{19}+2630\,{s}^{14}{t}^{18} \\
\end{eqnarray*}

\begin{eqnarray*}
&& +3315\,{s}^{15}{t}^{17}+3597\,{s}^{16}{t}^{16}+{s}^{9}{t}^{24}
+25\,{s}^{10}{t}^{23}+85\,{s}^{11}{t}^{22}+113\,{s}^{12}{t}^{21} \\
&& -44\,{s}^{13}{t}^{20}-589\,{s}^{14}{t}^{19}-1536\,{s}^{15}{t}^{18}
-2424\,{s}^{16}{t}^{17}+3\,{s}^{10}{t}^{24}-43\,{s}^{11}{t}^{23} \\
&& -310\,{s}^{12}{t}^{22}-1085\,{s}^{13}{t}^{21}-2662\,{s}^{14}{t}^{20}
-4897\,{s}^{15}{t}^{19}-7265\,{s}^{16}{t}^{18} \\
&& -8700\,{s}^{17}{t}^{17}-{s}^{10}{t}^{25}-47\,{s}^{11}{t}^{24}
-323\,{s}^{12}{t}^{23}-1104\,{s}^{13}{t}^{22}-2730\,{s}^{14}{t}^{21} \\
&& -5139\,{s}^{15}{t}^{20}-7802\,{s}^{16}{t}^{19}
-9930\,{s}^{17}{t}^{18}-22\,{s}^{11}{t}^{25}-158\,{s}^{12}{t}^{24} \\
&& -581\,{s}^{13}{t}^{23}-1465\,{s}^{14}{t}^{22}-2797\,{s}^{15}{t}^{21}
-4285\,{s}^{16}{t}^{20}-5471\,{s}^{17}{t}^{19} \\
&& -6163\,{s}^{18}{t}^{18}-5\,{s}^{11}{t}^{26}-39\,{s}^{12}{t}^{25}
-85\,{s}^{13}{t}^{24}-135\,{s}^{14}{t}^{23}+38\,{s}^{15}{t}^{22} \\
&& +557\,{s}^{16}{t}^{21}+1452\,{s}^{17}{t}^{20}+2361\,{s}^{18}{t}^{19}
-{s}^{11}{t}^{27}+5\,{s}^{12}{t}^{26}+103\,{s}^{13}{t}^{25} \\
&& +519\,{s}^{14}{t}^{24}+1588\,{s}^{15}{t}^{23}+3626\,{s}^{16}{t}^{22}
+6408\,{s}^{17}{t}^{21}+9221\,{s}^{18}{t}^{20} \\
&& +10970\,{s}^{19}{t}^{19}+6\,{s}^{12}{t}^{27}+97\,{s}^{13}{t}^{26}
+476\,{s}^{14}{t}^{25}+1579\,{s}^{15}{t}^{24}+3698\,{s}^{16}{t}^{23} \\
&& +6873\,{s}^{17}{t}^{22}+10313\,{s}^{18}{t}^{21}+12991\,{s}^{19}{t}^{20}
+2\,{s}^{12}{t}^{28}+38\,{s}^{13}{t}^{27}+213\,{s}^{14}{t}^{26} \\
&& +745\,{s}^{15}{t}^{25}+1914\,{s}^{16}{t}^{24}+3672\,{s}^{17}{t}^{23}
+5792\,{s}^{18}{t}^{22}+7543\,{s}^{19}{t}^{21} \\
&& +8636\,{s}^{20}{t}^{20}+8\,{s}^{13}{t}^{28}+26\,{s}^{14}{t}^{27}
+64\,{s}^{15}{t}^{26}+49\,{s}^{16}{t}^{25}-82\,{s}^{17}{t}^{24} \\
&& -581\,{s}^{18}{t}^{23}-1285\,{s}^{19}{t}^{22}-1988\,{s}^{20}{t}^{21}
-22\,{s}^{14}{t}^{28}-179\,{s}^{15}{t}^{27}-735\,{s}^{16}{t}^{26} \\
&& -2110\,{s}^{17}{t}^{25}-4476\,{s}^{18}{t}^{24}-7718\,{s}^{19}{t}^{23}
-10727\,{s}^{20}{t}^{22}-12706\,{s}^{21}{t}^{21} \\
&& -16\,{s}^{14}{t}^{29}-133\,{s}^{15}{t}^{28}-624\,{s}^{16}{t}^{27}
-1950\,{s}^{17}{t}^{26}-4569\,{s}^{18}{t}^{25}-8302\,{s}^{19}{t}^{24} \\
&& -12459\,{s}^{20}{t}^{23}-15619\,{s}^{21}{t}^{22}-4\,{s}^{14}{t}^{30}
-45\,{s}^{15}{t}^{29}-239\,{s}^{16}{t}^{28}-862\,{s}^{17}{t}^{27} \\
&& -2228\,{s}^{18}{t}^{26}-4481\,{s}^{19}{t}^{25}-7120\,{s}^{20}{t}^{24}
-9652\,{s}^{21}{t}^{23}-11078\,{s}^{22}{t}^{22} \\
&& -3\,{s}^{15}{t}^{30}-9\,{s}^{16}{t}^{29}-10\,{s}^{17}{t}^{28}
+4\,{s}^{18}{t}^{27}+100\,{s}^{19}{t}^{26}+291\,{s}^{20}{t}^{25} \\
&& +665\,{s}^{21}{t}^{24}+850\,{s}^{22}{t}^{23}+3\,{s}^{15}{t}^{31}
+41\,{s}^{16}{t}^{30}+238\,{s}^{17}{t}^{29}+892\,{s}^{18}{t}^{28} \\
&& +2388\,{s}^{19}{t}^{27}+4918\,{s}^{20}{t}^{26}+8085\,{s}^{21}{t}^{25}
+11130\,{s}^{22}{t}^{24}+12956\,{s}^{23}{t}^{23}.
\end{eqnarray*}

The numerator $N(T_{6,2};s,t)$ is a polynomial in $s$ and $t$ of
total degree 90, with leading term $(st)^{45}$,
having in total 854 terms. By using its symmetry properties,
we can write it as follows
\begin{eqnarray*}
&& N(T_{6,2};s,t)=(1+st)\left[ g_1(st)+(st)^{44} g_1(s^{-1}t^{-1}) \right] \\
&& \quad\quad -st^2\left[ g_2(s,t)+(st)^{42}g_2(t^{-1},s^{-1}) \right]
-s^2t\left[ g_2(t,s)+(st)^{42}g_2(s^{-1},t^{-1}) \right],
\end{eqnarray*}
where
\begin{eqnarray*}
g_1(x) &=& 1-2\,x+2\,{x}^{2}+4\,{x}^{3}-3\,{x}^{4}-6\,{x}^{5}+33\,{x}^{6}
+32\,{x}^{7}-59\,{x}^{8}-22\,{x}^{9} \\
&& +166\,{x}^{10}-8\,{x}^{11}-360\,{x}^{12}-110\,{x}^{13}+494\,{x}^{14}
+28\,{x}^{15}-711\,{x}^{16} \\
&& +192\,{x}^{17}+1203\,{x}^{18} +80\,{x}^{19}-1282\,{x}^{20}
-188\,{x}^{21}+528\,{x}^{22}
\end{eqnarray*}
and the polynomial $g_2$ has 206 terms:

\newpage
\begin{eqnarray*}
&& g_2(s,t)= \\
&& 1+t-3\,st-4\,s{t}^{2}-4\,s{t}^{3}-6\,{s}^{2}{t}^{2}-2\,s{t}^{4}-2\,{s}^{2}{t}^{3}
-s{t}^{5}+3\,{s}^{2}{t}^{4}+7\,{s}^{3}{t}^{3} \\
&& +4\,{s}^{2}{t}^{5}+9\,{s}^{3}{t}^{4}+4\,{s}^{2}{t}^{6}+5\,{s}^{3}{t}^{5}
+{s}^{4}{t}^{4}+2\,{s}^{2}{t}^{7}-3\,{s}^{3}{t}^{6}-19\,{s}^{4}{t}^{5}+{s}^{2}{t}^{8} \\
&& -6\,{s}^{3}{t}^{7}-30\,{s}^{4}{t}^{6}-52\,{s}^{5}{t}^{5}-7\,{s}^{3}{t}^{8}
-28\,{s}^{4}{t}^{7}-53\,{s}^{5}{t}^{6}-4\,{s}^{3}{t}^{9}-17\,{s}^{4}{t}^{8} \\
&& -27\,{s}^{5}{t}^{7}-35\,{s}^{6}{t}^{6}-2\,{s}^{3}{t}^{10}-3\,{s}^{4}{t}^{9}
+6\,{s}^{5}{t}^{8}+22\,{s}^{6}{t}^{7}+4\,{s}^{4}{t}^{10} +27\,{s}^{5}{t}^{9} \\
&& +56\,{s}^{6}{t}^{8}+78\,{s}^{7}{t}^{7}+6\,{s}^{4}{t}^{11}+29\,{s}^{5}{t}^{10}
+60\,{s}^{6}{t}^{9}+73\,{s}^{7}{t}^{8}+3\,{s}^{4}{t}^{12}+19\,{s}^{5}{t}^{11} \\
&& +31\,{s}^{6}{t}^{10}+19\,{s}^{7}{t}^{9}-5\,{s}^{8}{t}^{8}+{s}^{4}{t}^{13}
+7\,{s}^{5}{t}^{12}+{s}^{6}{t}^{11}-44\,{s}^{7}{t}^{10}-104\,{s}^{8}{t}^{9} \\
&& -{s}^{5}{t}^{13}-18\,{s}^{6}{t}^{12}-71\,{s}^{7}{t}^{11}-149\,{s}^{8}{t}^{10}
-213\,{s}^{9}{t}^{9}-{s}^{5}{t}^{14}-17\,{s}^{6}{t}^{13}-56\,{s}^{7}{t}^{12} \\
&& -102\,{s}^{8}{t}^{11}-148\,{s}^{9}{t}^{10}-{s}^{5}{t}^{15}
-8\,{s}^{6}{t}^{14}-20\,{s}^{7}{t}^{13}-5\,{s}^{8}{t}^{12}+29\,{s}^{9}{t}^{11} \\
&& +67\,{s}^{10}{t}^{10}-2\,{s}^{6}{t}^{15}+10\,{s}^{7}{t}^{14}+73\,{s}^{8}{t}^{13}
+191\,{s}^{9}{t}^{12}+318\,{s}^{10}{t}^{11}+{s}^{6}{t}^{16} \\
&& +18\,{s}^{7}{t}^{15}+100\,{s}^{8}{t}^{14}+246\,{s}^{9}{t}^{13}+428\,{s}^{10}{t}^{12}
+543\,{s}^{11}{t}^{11}+12\,{s}^{7}{t}^{16}+69\,{s}^{8}{t}^{15} \\
&& +185\,{s}^{9}{t}^{14}+324\,{s}^{10}{t}^{13}+434\,{s}^{11}{t}^{12}
+4\,{s}^{7}{t}^{17}+29\,{s}^{8}{t}^{16}+68\,{s}^{9}{t}^{15}+111\,{s}^{10}{t}^{14} \\
&& +109\,{s}^{11}{t}^{13}+96\,{s}^{12}{t}^{12}-20\,{s}^{9}{t}^{16}
-86\,{s}^{10}{t}^{15}-198\,{s}^{11}{t}^{14}-335\,{s}^{12}{t}^{13} \\
&& -6\,{s}^{8}{t}^{18}-49\,{s}^{9}{t}^{17}-155\,{s}^{10}{t}^{16}-331\,{s}^{11}{t}^{15}
-519\,{s}^{12}{t}^{14}-653\,{s}^{13}{t}^{13}-3\,{s}^{8}{t}^{19} \\
&& -34\,{s}^{9}{t}^{18}-118\,{s}^{10}{t}^{17}-247\,{s}^{11}{t}^{16}
-390\,{s}^{12}{t}^{15}-474\,{s}^{13}{t}^{14}-{s}^{8}{t}^{20}-13\,{s}^{9}{t}^{19} \\
&& -45\,{s}^{10}{t}^{18}-79\,{s}^{11}{t}^{17}-82\,{s}^{12}{t}^{16}
-30\,{s}^{13}{t}^{15}+30\,{s}^{14}{t}^{14}-2\,{s}^{9}{t}^{20}+4\,{s}^{10}{t}^{19} \\
&& +49\,{s}^{11}{t}^{18}+169\,{s}^{12}{t}^{17}+377\,{s}^{13}{t}^{16}
+575\,{s}^{14}{t}^{15}+11\,{s}^{10}{t}^{20}+74\,{s}^{11}{t}^{19}  \\
&& +220\,{s}^{12}{t}^{18}+469\,{s}^{13}{t}^{17}+727\,{s}^{14}{t}^{16}
+905\,{s}^{15}{t}^{15}+6\,{s}^{10}{t}^{21}+39\,{s}^{11}{t}^{20} \\
&& +124\,{s}^{12}{t}^{19}+264\,{s}^{13}{t}^{18}+405\,{s}^{14}{t}^{17}
+503\,{s}^{15}{t}^{16}+{s}^{10}{t}^{22}+6\,{s}^{11}{t}^{21}-{s}^{12}{t}^{20} \\
&& -26\,{s}^{13}{t}^{19}-124\,{s}^{14}{t}^{18}-262\,{s}^{15}{t}^{17}
-399\,{s}^{16}{t}^{16}-5\,{s}^{11}{t}^{22}-51\,{s}^{12}{t}^{21} \\
&& -193\,{s}^{13}{t}^{20}-478\,{s}^{14}{t}^{19}-874\,{s}^{15}{t}^{18}
-1242\,{s}^{16}{t}^{17}-3\,{s}^{11}{t}^{23}-40\,{s}^{12}{t}^{22} \\
&& -174\,{s}^{13}{t}^{21}-491\,{s}^{14}{t}^{20}-954\,{s}^{15}{t}^{19}
-1452\,{s}^{16}{t}^{18}-1750\,{s}^{17}{t}^{17}-{s}^{11}{t}^{24} \\
&& -16\,{s}^{12}{t}^{23}-81\,{s}^{13}{t}^{22}-263\,{s}^{14}{t}^{21}
-578\,{s}^{15}{t}^{20}-935\,{s}^{16}{t}^{19}-1205\,{s}^{17}{t}^{18} \\
&& -3\,{s}^{12}{t}^{24}-10\,{s}^{13}{t}^{23}-33\,{s}^{14}{t}^{22}
-74\,{s}^{15}{t}^{21}-118\,{s}^{16}{t}^{20}-106\,{s}^{17}{t}^{19} \\
&& -77\,{s}^{18}{t}^{18}+10\,{s}^{13}{t}^{24}+59\,{s}^{14}{t}^{23}
+201\,{s}^{15}{t}^{22}+459\,{s}^{16}{t}^{21}+793\,{s}^{17}{t}^{20} \\
&& +1102\,{s}^{18}{t}^{19}+7\,{s}^{13}{t}^{25}+52\,{s}^{14}{t}^{24}
+205\,{s}^{15}{t}^{23}+542\,{s}^{16}{t}^{22}+1030\,{s}^{17}{t}^{21} \\
&& +1515\,{s}^{18}{t}^{20}+1801\,{s}^{19}{t}^{19}+2\,{s}^{13}{t}^{26}
+20\,{s}^{14}{t}^{25}+97\,{s}^{15}{t}^{24}+296\,{s}^{16}{t}^{23} \\
&& +647\,{s}^{17}{t}^{22}+1055\,{s}^{18}{t}^{21}+1346\,{s}^{19}{t}^{20}
+3\,{s}^{14}{t}^{26}+15\,{s}^{15}{t}^{25}+52\,{s}^{16}{t}^{24} \\
&& +113\,{s}^{17}{t}^{23}+208\,{s}^{18}{t}^{22}+272\,{s}^{19}{t}^{21}
+319\,{s}^{20}{t}^{20}-{s}^{14}{t}^{27}-8\,{s}^{15}{t}^{26} \\
&& -45\,{s}^{16}{t}^{25}-157\,{s}^{17}{t}^{24}-370\,{s}^{18}{t}^{23}
-621\,{s}^{19}{t}^{22}-799\,{s}^{20}{t}^{21}-3\,{s}^{15}{t}^{27} \\
&& -18\,{s}^{16}{t}^{26}-76\,{s}^{17}{t}^{25}-205\,{s}^{18}{t}^{24}
-403\,{s}^{19}{t}^{23}-572\,{s}^{20}{t}^{22}-672\,{s}^{21}{t}^{21}.
\end{eqnarray*}

\section{ Appendix B: Formanek's formulae }

Let $d$ be any positive integer.
Denote by $\mu$ a partition of a positive integer $k$ and by
$\chi^\mu$ the corresponding irreducible complex character of
the symmetric group $S_k$.
Define the {\em length}, $l(\mu)$, of $\mu$ to be
the number of parts of $\mu$.

Let $\Lambda_d$ denote the ring of symmetric polynomials in $d$ variables,
$t_1,\ldots,t_d$.
If $\mu$ has at most $d$ parts, we denote by $f_{\mu,d}\in\Lambda_d$
the corresponding Schur function.

Define the Frobenius map ${\rm Fr}_d$
to be the additive homomorphism from the direct sum of the character
rings of all $S_k$'s to $\Lambda_d$ by setting
${\rm Fr}_d(\chi^\mu)=f_{\mu,d}$ for each partition $\mu$ of $k$
having at most $d$ parts and ${\rm Fr}_d(\chi^\mu)=0$ otherwise.

We can now state Formanek's formulae for the algebras $C_{n,d}$
and $T_{n,d}$ (see \cite{BS})):
\begin{eqnarray}
P(C_{n,d};t_1,\ldots,t_d) &=& \label{form-1}
\sum_{k\ge0} {\rm Fr}_d\left( \theta_n^{(k)} \right)(t_1,\ldots,t_d), \\
P(T_{n,d};t_1,\ldots,t_d) &=& \label{form-2}
\sum_{k\ge0} {\rm Fr}_d\left( \theta_n^{(k+1)}\downarrow S_k \right)
(t_1,\ldots,t_d),
\end{eqnarray}
where $\theta_n^{(k)}$ is a particular character of $S_k$. This
character is defined by the formula
\begin{equation} \label{karakter}
\theta_n^{(k)} = \sum_{\mu:\, l(\mu)\le n} \chi^\mu\otimes\chi^\mu,
\end{equation}
where $\otimes$ is the usual tensor product of characters of $S_k$.

Recall that we are mainly interested in the case $d=2$. In that case
we set $s=t_1$, $t=t_2$, and $f_\mu=f_{\mu,2}$.
If $\mu=(p,q)$ with $p\ge q\ge0$, then
\begin{equation} \label{schur}
 f_\mu=f_{{p,q}}=(st)^{q}\left( s^{p-q}+s^{p-q-1}t+\cdots+st^{p-q-1}+t^{p-q}
\right).
\end{equation}
If $q=0$ we shall write $f_p=f_{{p,0}}$.

By using GAP \cite{GAP4}, we computed the first 26 terms of the
series (\ref{form-1}) in
the case $(n,d)=(5,2)$ and we obtained the formula displayed below.
By substituting the expressions (\ref{schur}) for the Schur functions
$f_\mu$ into this formula,
one obtains an initial chunk of the bivariate Taylor series
of $P(C_{5,2};s,t)$ which agrees with the bivariate Taylor series
of the rational function that we have computed
by using the Molien--Weyl formula (see Theorem \ref{Prva}).

The $k$-th summand in (\ref{form-1}), when written as a linear
combination of Schur functions, gives the decomposition of the
character of the representation of $\GL_d$ on the $k$-th
homogeneous component of $C_{n,d}$.

Let $c_{n,2}^\bullet(k)$ denote the sum of the coefficients of the
Schur functions $f_{k-p,p}$, $p=0,1,\ldots$, $2p\le k$,
in  (\ref{form-1}) when $d=2$.
This is the length of the homogeneous component of
$C_{n,2}$ of degree $k$ as a $\GL_2$-module.
Since $C_{n,2}^\bullet$ is the algebra of $\GL_1$-invariants
in $C_{n,2}$, the coefficient of $t^{2k}$ in the Taylor
series of $P(C_{n,2}^\bullet;t)$ must be equal to 
$c_{n,2}^\bullet(2k)$.
It is easy to verify this assertion when $n=5$ and $k\le12$
by using the expansion for $P(C_{5,2})$ below and the one
for $P(C_{5,2}^\bullet)$ at the end of Section \ref{manja-gr}.

\newpage
\begin{eqnarray*}
&& P(C_{5,2};s,t)= \\
&& 1+f_{{1}}+2\,f_{{2}}+f_{{2,1}}+3\,f_{{3}}+3\,f_{{2,2}}
+2\,f_{{3,1}}+5\,f_{{4}}+6\,f_{{3,2}}+5\,f_{{4,1}}
+7\,f_{{5}} \\
&& +4\,f_{{3,3}}+15\,f_{{4,2}}+8\,f_{{5,1}}+10\,f_{{6}}
+18\,f_{{4,3}}
+25\,f_{{5,2}}+14\,f_{{6,1}}+13\,f_{{7}}+24\,f_{{4,4}} \\
&& +37\,f_{{5,3}}+44\,f_{{6,2}}+20\,f_{{7,1}}+18\,f_{{8}}
+58\,f_{{5,4}}+76\,f_{{6,3}}+66\,f_{{7,2}}+30\,f_{{8,1}} \\
&& +23\,f_{{9}}+44\,f_{{5,5}}+136\,f_{{6,4}}+126\,f_{{7,3}}
+101\,f_{{8,2}}+41\,f_{{9,1}}+30\,f_{{10}}+163\,f_{{6,5}} \\
&& +246\,f_{{7,4}}+207\,f_{{8,3}}+142\,f_{{9,2}}+57\,f_{{10,1}}
+37\,f_{{11}}+171\,f_{{6,6}}+354\,f_{{7,5}} \\
&& +431\,f_{{8,4}}+311\,f_{{9,3}}+200\,f_{{10,2}}+74\,f_{{11,1}}
+47\,f_{{12}}+476\,f_{{7,6}}+700\,f_{{8,5}} \\
&& +681\,f_{{9,4}}+458\,f_{{10,3}}+267\,f_{{11,2}}+98\,f_{{12,1}}
+57\,f_{{13}}+388\,f_{{7,7}}+1080\,f_{{8,6}} \\
&& +1204\,f_{{9,5}}+1047\,f_{{10,4}}+640\,f_{{11,3}}+357\,f_{{12,2}}
+124\,f_{{13,1}}+70\,f_{{14}}+1277\,f_{{8,7}} \\
&& +2024\,f_{{9,6}}+1973\,f_{{10,5}}+1517\,f_{{11,4}}
+884\,f_{{12,3}}+460\,f_{{13,2}}+157\,f_{{14,1}}+84\,f_{{15}} \\
&& +1166\,f_{{8,8}}+2811\,f_{{9,7}}+3534\,f_{{10,6}}
+3014\,f_{{11,5}}+2160\,f_{{12,4}}+1177\,f_{{13,3}} \\
&& +591\,f_{{14,2}}+194\,f_{{15,1}}+101\,f_{{16}}+3419\,f_{{9,8}}
+5427\,f_{{10,7}}+5682\,f_{{11,6}} \\
&& +4470\,f_{{12,5}}+2957\,f_{{13,4}}+1550\,f_{{14,3}}
+740\,f_{{15,2}}+240\,f_{{16,1}}+119\,f_{{17}} \\
&& +2808\,f_{{9,9}}+7592\,f_{{10,8}}+9371\,f_{{11,7}}
+8780\,f_{{12,6}}+6352\,f_{{13,5}}+3985\,f_{{14,4}} \\
&& +1992\,f_{{15,3}}+924\,f_{{16,2}}+290\,f_{{17,1}}
+141\,f_{{18}}+8693\,f_{{10,9}}+14284\,f_{{11,8}} \\
&& +15284\,f_{{12,7}}+12921\,f_{{13,6}}+8823\,f_{{14,5}}
+5232\,f_{{15,4}}+2535\,f_{{16,3}} \\
&& +1131\,f_{{17,2}}+351\,f_{{18,1}}+164\,f_{{19}}
+7344\,f_{{10,10}}+18927\,f_{{11,9}}+24781\,f_{{12,8}} \\
&& +23520\,f_{{13,7}}+18489\,f_{{14,6}}+11909\,f_{{15,5}}
+6784\,f_{{16,4}}+3167\,f_{{17,3}}+1380\,f_{{18,2}} \\
&& +417\,f_{{19,1}}+192\,f_{{20}}+21565\,f_{{11,10}}
+35929\,f_{{12,9}}+40009\,f_{{13,8}}+34897\,f_{{14,7}} \\
&& +25625\,f_{{15,6}}+15798\,f_{{16,5}}+8622\,f_{{17,4}}
+3926\,f_{{18,3}}+1658\,f_{{19,2}}+496\,f_{{20,1}} \\
&& +221\,f_{{21}}+17281\,f_{{11,11}}+46991\,f_{{12,10}}
+61801\,f_{{13,9}}+61722\,f_{{14,8}} \\
&& +49917\,f_{{15,7}}+34778\,f_{{16,6}}+20520\,f_{{17,5}}
+10849\,f_{{18,4}}+4796\,f_{{19,3}}+1986\,f_{{20,2}} \\
&& +582\,f_{{21,1}}+255\,f_{{22}}+51694\,f_{{12,11}}
+87853\,f_{{13,10}}+100058\,f_{{14,9}}+91235\,f_{{15,8}} \\
&& +69582\,f_{{16,7}}+46117\,f_{{17,6}}+26294\,f_{{18,5}}
+13444\,f_{{19,4}}+5820\,f_{{20,3}}+2350\,f_{{21,2}} \\
&& +682\,f_{{22,1}}+291\,f_{{23}}+41569\,f_{{12,12}}
+111058\,f_{{13,11}}+150865\,f_{{14,10}} \\
&& +153818\,f_{{15,9}}+130796\,f_{{16,8}}+94507\,f_{{17,7}}
+60179\,f_{{18,6}}+33154\,f_{{19,5}} \\
&& +16519\,f_{{20,4}}+6983\,f_{{21,3}}+2772\,f_{{22,2}}
+790\,f_{{23,1}}+333\,f_{{24}}+120672\,f_{{13,12}} \\
&& +207439\,f_{{14,11}}+242629\,f_{{15,10}}+227776\,f_{{16,9}}
+182080\,f_{{17,8}}+125907\,f_{{18,7}} \\
&& +77153\,f_{{19,6}}+41349\,f_{{20,5}}+20055\,f_{{21,4}}
+8328\,f_{{22,3}}+3237\,f_{{23,2}} \\
&& +915\,f_{{24,1}}+377\,f_{{25}}+\cdots.
\end{eqnarray*}

\section{ Appendix C: Taylor expansions of $P(C_{n,2};t)$
and $P(C_{n,2}(0);t)$} 

In Table 9 we give the Taylor series of
$P(C_{n,2};t)$ and $P(C_{n,2}(0);t)$ for $n\le6$
including the terms of degree $k<20$.
In the former case we also tabulate (see Table 10) the coefficients
for $k\le12$ and make a couple of observations. We are then lead
to some speculations concerning certain limits of
the algebras $C_{n,d}$, which we are going to introduce now.

\begin{center}

{\bf Table 9: Taylor expansions of $P(C_{n,2};t)$ and $P(C_{n,2}(0);t)$}
\begin{eqnarray*}
P(C_{0,2};t) &=& 1, \\
P(C_{1,2};t) &=&
1+2\,t+3\,{t}^{2}+4\,{t}^{3}+5\,{t}^{4}+6\,{t}^{5}+7\,{t}^{6}+8\,
{t}^{7}+9\,{t}^{8}+10\,{t}^{9} \\
&& +11\,{t}^{10}+12\,{t}^{11}+13\,{t}^{12}
+14\,{t}^{13}+15\,{t}^{14}+16\,{t}^{15}+17\,{t}^{16} \\
&& +18\,{t}^{17}+19\,{t}^{18}+20\,{t}^{19}+\cdots, \\
P(C_{1,2}(0);t) &=& 1, \\
P(C_{2,2};t) &=&
1+2\,t+6\,{t}^{2}+10\,{t}^{3}+20\,{t}^{4}+30\,{t}^{5}+50\,{t}^{6}+70\,
{t}^{7}+105\,{t}^{8} \\
&& +140\,{t}^{9}+196\,{t}^{10}+252\,{t}^{11}+336\,{t}
^{12}+420\,{t}^{13}+540\,{t}^{14}+660\,{t}^{15} \\
&& +825\,{t}^{16}+990\,{t}^{17}+1210\,{t}^{18}+1430\,{t}^{19}+\cdots, \\
P(C_{2,2}(0);t) &=&
1+3\,{t}^{2}+6\,{t}^{4}+10\,{t}^{6}+15\,{t}^{8}+21\,{t}^{10}+28\,{t}^{12}
+36\,{t}^{14}+45\,{t}^{16} \\
&& +55\,{t}^{18}+\cdots, \\
P(C_{3,2};t) &=&
1+2\,t+6\,{t}^{2}+14\,{t}^{3}+29\,{t}^{4}+56\,{t}^{5}+107\,{t}^{6}+186
\,{t}^{7}+320\,{t}^{8} \\
&& +530\,{t}^{9}+851\,{t}^{10}+1332\,{t}^{11}+2051\,{t}^{12}+3074\,{t}^{13}
+4544\,{t}^{14} \\
&& +6602\,{t}^{15}+9444\,{t}^{16}+13322\,{t}^{17}+18579\,{t}^{18}
+25564\,{t}^{19}+\cdots, \\
P(C_{3,2}(0);t) &=&
1+3\,{t}^{2}+4\,{t}^{3}+7\,{t}^{4}+12\,{t}^{5}+24\,{t}^{6}+28\,{t}^{7}
+55\,{t}^{8}+76\,{t}^{9} \\
&& +111\,{t}^{10}+160\,{t}^{11}+238\,{t}^{12}+304
\,{t}^{13}+447\,{t}^{14}+588\,{t}^{15} \\
&& +784\,{t}^{16}+1036\,{t}^{17}+1379\,{t}^{18}+1728\,{t}^{19}+\cdots, \\
P(C_{4,2};t) &=&
1+2\,t+6\,{t}^{2}+14\,{t}^{3}+34\,{t}^{4}+68\,{t}^{5}+144\,{t}^{6}+276
\,{t}^{7}+534\,{t}^{8} \\
&& +974\,{t}^{9}+1774\,{t}^{10}+3106\,{t}^{11}+5410\,{t}^{12}
+9146\,{t}^{13}+15334\,{t}^{14} \\
&& +25158\,{t}^{15}+40884\,{t}^{16}+65264\,{t}^{17}+103204\,{t}^{18} \\
&& +160808\,{t}^{19}+\cdots, \\
P(C_{4,2}(0);t) &=&
1+3\,{t}^{2}+4\,{t}^{3}+12\,{t}^{4}+14\,{t}^{5}+42\,{t}^{6}+56\,{t}^{7
}+126\,{t}^{8}+182\,{t}^{9} \\
&& +360\,{t}^{10}+532\,{t}^{11}+972\,{t}^{12}+
1432\,{t}^{13}+2452\,{t}^{14}+3636\,{t}^{15} \\
&& +5902\,{t}^{16}+8654\,{t}^{17}+13560\,{t}^{18}+19664\,{t}^{19}+\cdots, \\
P(C_{5,2};t) &=&
1+2\,t+6\,{t}^{2}+14\,{t}^{3}+34\,{t}^{4}+74\,{t}^{5}+159\,{t}^{6}+324
\,{t}^{7}+657\,{t}^{8} \\
&& +1286\,{t}^{9}+2488\,{t}^{10}+4702\,{t}^{11}+
8790\,{t}^{12}+16146\,{t}^{13} \\
&& +29326\,{t}^{14}+52526\,{t}^{15}+93064\,{t}^{16}+162910\,{t}^{17}
+282267\,{t}^{18} \\
&& +483792\,{t}^{19}+\cdots, \\
P(C_{5,2}(0);t) &=&
1+3\,{t}^{2}+4\,{t}^{3}+12\,{t}^{4}+20\,{t}^{5}+45\,{t}^{6}+80\,{t}^{7
}+168\,{t}^{8}+296\,{t}^{9} \\
&& +573\,{t}^{10}+1012\,{t}^{11}+1874\,{t}^{12
}+3268\,{t}^{13}+5824\,{t}^{14}+10020\,{t}^{15} \\
&& +17338\,{t}^{16}+29308\,{t}^{17}+49511\,{t}^{18}+82168\,{t}^{19}+\cdots, \\
P(C_{6,2};t) &=&
1+2\,t+6\,{t}^{2}+14\,{t}^{3}+34\,{t}^{4}+74\,{t}^{5}+166\,{t}^{6}+342
\,{t}^{7}+716\,{t}^{8} \\
&& +1442\,{t}^{9}+2898\,{t}^{10}+5686\,{t}^{11}+
11122\,{t}^{12}+21366\,{t}^{13} \\
&& +40842\,{t}^{14}+77098\,{t}^{15}+144581\,{t}^{16}+268376\,{t}^{17}
+494812\,{t}^{18} \\
&& +904056\,{t}^{19}+\cdots, \\
P(C_{6,2}(0);t) &=&
1+3\,{t}^{2}+4\,{t}^{3}+12\,{t}^{4}+20\,{t}^{5}+52\,{t}^{6}+84\,{t}^{7
}+198\,{t}^{8}+352\,{t}^{9} \\
&& +730\,{t}^{10}+1332\,{t}^{11}+2648\,{t}^{12
}+4808\,{t}^{13}+9232\,{t}^{14}+16780\,{t}^{15} \\
&& +31227\,{t}^{16}+56312\,{t}^{17}+102641\,{t}^{18}+182808\,{t}^{19}+\cdots.
\end{eqnarray*}

\end{center}

Define the $\bZ^d$-graded algebra $C_{\infty,d}$ as the inverse limit of
\[ C_{1,d} \leftarrow C_{2,d} \leftarrow C_{3,d} \leftarrow \, \cdots. \]
By adapting a definition of Formanek \cite[p. 52]{EF2},
we refer to $C_{\infty,d}$
as the {\em pure free trace ring} on $d$ generators.
One can next take the direct limit of
\[ C_{\infty,1} \to C_{\infty,2} \to C_{\infty,3} \to \cdots \]
to obtain the $\bZ^\infty$-graded algebra $C_{\infty,\infty}$,
which is the pure free trace ring
on countably many generators $x_1,x_2,x_3,\ldots$.
It follows from the Second Fundamental Theorem for invariants
of $n\times n$ matrices (see e.g. \cite[Theorem 50]{EF2})
that $C_{\infty,\infty}$ is indeed isomorphic to the pure
free trace ring as defined by Formanek.

Let us write $c_{n,2}(k)$ for the coefficients of the Poincar\'{e} series
\[ P(C_{n,2};t)=\sum_{k\ge0} c_{n,2}(k) t^k \]
and let us display these coefficients in an infinite table
(rows indexed by $n\ge0$ and columns by $k\ge0$). The data from
Table 9 give the top portion of Table 10.

\begin{center}

{\bf Table 10: The coefficients $c_{n,2}(k)$}
\[
\begin{array}{rrrrrrrrrrrrrr}
1 & 0 & 0 & 0 & 0 & 0 & 0 & 0 & 0 & 0 & 0 & 0 & 0 & \cdots \\
1 & 2 & 3 & 4 & 5 & 6 & 7 & 8 & 9 & 10 & 11 & 12 & 13 & \\
1 & 2 & 6 & 10 & 20 & 30 & 50 & 70 & 105 & 140 & 196 & 252 & 336 & \\
1 & 2 & 6 & 14 & 29 & 56 & 107 & 186 & 320 & 530 & 851 & 1332 & 2051 & \\
1 & 2 & 6 & 14 & 34 & 68 & 144 & 276 & 534 & 974 & 1774 & 3106 & 5410 & \\
1 & 2 & 6 & 14 & 34 & 74 & 159 & 324 & 657 & 1286 & 2488 & 4702 & 8790 & \\
1 & 2 & 6 & 14 & 34 & 74 & 166 & 342 & 716 & 1442 & 2898 & 5686 & 11122 & \\
\hline
1 & 2 & 6 & 14 & 34 & 74 & 166 & 350 & 737 & 1512 & 3087 & 6194 & 12376 & \\
1 & 2 & 6 & 14 & 34 & 74 & 166 & 350 & 746 & 1536 & 3168 & 6416 & 12982 & \\
1 & 2 & 6 & 14 & 34 & 74 & 166 & 350 & 746 & 1546 & 3195 & 6508 & 13237 & \\
1 & 2 & 6 & 14 & 34 & 74 & 166 & 350 & 746 & 1546 & 3206 & 6538 & 13340 & \\
1 & 2 & 6 & 14 & 34 & 74 & 166 & 350 & 746 & 1546 & 3206 & 6550 & 13373 & \\
1 & 2 & 6 & 14 & 34 & 74 & 166 & 350 & 746 & 1546 & 3206 & 6550 & 13386 & \\
\vdots &&&&&&&&&&&&& \end{array}
\]

\end{center}

We observed from the top part of this table that apparently each column
stabilizes and that
\[ \lim_{n\to\infty} P(C_{n,2};t)=\sum_{k\ge0} c_{n,2}(n)t^n . \]
After making this observation, we looked up the diagonal sequence
\[ \{ c_{n,2}(n) \}_{n\ge0} = 1, 2, 6, 14, 34, 74, 166, \ldots \]
in the On-Line Encyclopedia of Integer Sequences \cite{NS}
by entering only these 7 integers. We were pleasantly surprised
to find that it is registered there as the sequence A070933, and
identified as the sequence of coefficients in the power series
expansion of the infinite product
\[ \prod_{k\ge1} \frac{1}{1-2t^k}. \]
The first 30 terms of A070933 are listed in \cite{NS}.
The above infinite product should be the Poincar\'{e} series of 
the algebra $C_{\infty,2}$ (see the definition below). 
In the bigraded case one would have to replace the above product with
\[ \prod_{k\ge1} \frac{1}{1-s^k-t^k}. \]

More generally, we expect that the multigraded Poincar\'{e} series of
$C_{\infty,d}$ and $C_{\infty,\infty}$ be given by
\[ P(C_{\infty,d};t_1,\ldots,t_d)=\prod_{k\ge1}
\frac{1}{1-t_1^k-\cdots-t_d^k} \]
and
\begin{equation} \label{besk-pr}
 P(C_{\infty,\infty};t_1,t_2,\ldots)=\prod_{k\ge1}
\frac{1}{1-p_k},
\end{equation}
respectively, where the $p_k$ are the usual power sum symmetric functions:
\[ p_k=t_1^k+t_2^k+\cdots. \]

The latter formula is indeed valid. As explained in \cite{EF1},
it follows from the Procesi--Razmyslov theorem that
\[ P(C_{n,\infty};t)=\sum_{\mu;\, l(\mu)\le n}
 {\rm Fr}(\chi^\mu \otimes \chi^\mu). \]
The Frobenius map Fr is the additive map from the direct sum
of the character rings of all $S_k$'s to the ring, $\Lambda$,
of symmetric functions in infinitely many variables $t_1,t_2,\ldots$.
It is defined by setting, for all partitions $\mu$,
${\rm Fr}(\chi^\mu)=f_\mu\in\Lambda$, the Schur function
corresponding to the partition $\mu$.
By letting $n\to\infty$, we obtain that
\[ P(C_{\infty,\infty};t)=\sum_\mu {\rm Fr}(\chi^\mu \otimes \chi^\mu), \]
where the summation is now over all partitions $\mu$.
It remains to observe that the right hand side of this formula and
the one of (\ref{besk-pr}) are equal, see Macdonald's classic
\cite[Chapter 1, Section 7, Example 9(a)]{IM}.

Another interesting observation about the above table is that
apparently the second differences
\begin{equation*}
\alpha_k = c_{n,2}(n+k)-c_{n-1,2}(n+k)-c_{n-1,2}(n+k-1)+c_{n-2,2}(n+k-1)  
\end{equation*}
are independent of $n$ for $n\ge k$. The sequence
\begin{equation} \label{myst}
 \{ \alpha_{k} \}_{k\ge0} = 1,\, 3,\, 11,\, 33,\, 98,\, 270,\, \ldots
\end{equation}
has not been recorded so far in \cite{NS}.
By using the sequence A070933 and the hypothetical rules mentioned above,
we extended the top portion of Table 10 with 6 additional rows.
(Actually the above rules allow us to fill in the next three columns of
the table, but the page is not wide enough.)
Subsequently, we verified the validity of these additional rows by using
Formanek's formula, and at the same time we enlarged the number of
columns to 26, i.e., $0\le k\le25$. This made it possible to compute
a few more terms of the sequence (\ref{myst}):
\[ 1,\, 3,\, 11,\, 33,\, 98,\, 270,\, 736,\, 1932,\, 5009,\, 
12727,\, 31977,\, 79307,\, 194947,\,  \ldots \]

This sequence has another hypothetical incarnation. Recall the
integers $c_{n,2}^\bullet(k)$ introduced in Appendix B.
We conjecture that, for each fixed $k\ge0$,
\[  \alpha_k=c_{n,2}^\bullet(n+k)-c_{n-1,2}^\bullet(n+k) \]
is valid for sufficiently large $n$.

\section{ Appendix D: Taylor expansion of $P(T_{n,2};t)$} 

In this appendix we give, for $n\le6$, the coefficients of $t^k$, $k<20$,
in the power series expansions of $P(T_{n,2};t)$.
We also tabulate the coefficients for $k\le12$
and make some interesting observations about the table.

\begin{center}

{\bf Table 11: Taylor expansions of $P(T_{n,2};t)$}
\begin{eqnarray*}
P(T_{0,2};t) &=& 1, \\
P(T_{1,2};t) &=& 1+2\,t+3\,{t}^{2}+4\,{t}^{3}+5\,{t}^{4}+6\,{t}^{5}+7\,{t}^{6}
+8\,{t}^{7}+9\,{t}^{8}+10\,{t}^{9} \\
&& +11\,{t}^{10}+12\,{t}^{11}+13\,{t}^{12}+14\,{t}^{13}
+15\,{t}^{14}+16\,{t}^{15}+17\,{t}^{16} \\
&& +18\,{t}^{17}+19\,{t}^{18}+20\,{t}^{19} +\cdots, \\
P(T_{2,2};t) &=& 1+4\,t+11\,{t}^{2}+24\,{t}^{3}+46\,{t}^{4}+80\,{t}^{5}
+130\,{t}^{6}+200\,{t}^{7}+295\,{t}^{8} \\
&& +420\,{t}^{9}+581\,{t}^{10}+784\,{t}^{11}+1036\,{t}^{12}
+1344\,{t}^{13}+1716\,{t}^{14} \\
&& +2160\,{t}^{15}+2685\,{t}^{16}+3300\,{t}^{17}+4015\,{t}^{18}
+4840\,{t}^{19} +\cdots, \\
P(T_{3,2};t) &=& 1+4\,t+14\,{t}^{2}+38\,{t}^{3}+93\,{t}^{4}+204\,{t}^{5}
+419\,{t}^{6}+806\,{t}^{7}+1480\,{t}^{8} \\
&& +2600\,{t}^{9}+4411\,{t}^{10}+7244\,{t}^{11}
+11579\,{t}^{12}+18048\,{t}^{13}+27528\,{t}^{14} \\
&& +41150\,{t}^{15}+60428\,{t}^{16}+87280\,{t}^{17}
+124203\,{t}^{18}+174308\,{t}^{19} +\cdots, \\
P(T_{4,2};t) &=& 1+4\,t+14\,{t}^{2}+42\,{t}^{3}+113\,{t}^{4}+278\,{t}^{5}
+646\,{t}^{6}+1418\,{t}^{7}+2979\,{t}^{8} \\
&& +6018\,{t}^{9}+11752\,{t}^{10}+22256\,{t}^{11}
+41030\,{t}^{12}+73784\,{t}^{13} \\
&& +129748\,{t}^{14}+223498\,{t}^{15}+377753\,{t}^{16}+627314\,{t}^{17} \\
&& +1024882\,{t}^{18}+1649026\,{t}^{19} +\cdots, \\
P(T_{5,2};t) &=& 1+4\,t+14\,{t}^{2}+42\,{t}^{3}+118\,{t}^{4}
+304\,{t}^{5}+747\,{t}^{6}+1748\,{t}^{7}+3949\,{t}^{8} \\
&& +8620\,{t}^{9}+18296\,{t}^{10}+37818\,{t}^{11}+76398\,{t}^{12}
+151022\,{t}^{13} \\
&& +292754\,{t}^{14}+557130\,{t}^{15}+1042364\,{t}^{16}+1919044\,{t}^{17} \\
&& +3480203\,{t}^{18}+6221668\,{t}^{19} +\cdots, \\
P(T_{6,2};t) &=& 1+4\,t+14\,{t}^{2}+42\,{t}^{3}+118\,{t}^{4}+310\,{t}^{5}
+779\,{t}^{6}+1876\,{t}^{7}+4382\,{t}^{8} \\
&& +9948\,{t}^{9}+22057\,{t}^{10}+47850\,{t}^{11}
+101844\,{t}^{12}+212946\,{t}^{13} \\
&& +438118\,{t}^{14}+887814\,{t}^{15}+1773827\,{t}^{16}+3496850\,{t}^{17} \\
&& +6806682\,{t}^{18}+13089748\,{t}^{19} +\cdots.
\end{eqnarray*}

\end{center}

Let us write $d_{n,2}(k)$ for the coefficients of the Poincar\'{e} series
\[ P(T_{n,2};t)=\sum_{k\ge0} d_{n,2}(k) t^k \]
and let us display these coefficients in an infinite table
(rows indexed by $n\ge0$ and columns by $k\ge0$).

\begin{center}

{\bf Table 12: The coefficients $d_{n,2}(k)$}
\[
\begin{array}{rrrrrrrrrrrrrr}
1 & 0 & 0 & 0 & 0 & 0 & 0 & 0 & 0 & 0 & 0 & 0 & 0 & \cdots \\
1 & 2 & 3 & 4 & 5 & 6 & 7 & 8 & 9 & 10 & 11 & 12 & 13 & \\
1 & 4 & 11 & 24 & 46 & 80 & 130 & 200 & 295 & 420 & 581 & 784 & 1036 & \\
1 & 4 & 14 & 38 & 93 & 204 & 419 & 806 & 1480 & 2600 & 4411 & 7244 & 11579 & \\
1 & 4 & 14 & 42 & 113 & 278 & 646 & 1418 & 2979 & 6018 & 11752 & 22256 & 41030 & \\
1 & 4 & 14 & 42 & 118 & 304 & 747 & 1748 & 3949 & 8620 & 18296 & 37818 & 76398 & \\
1 & 4 & 14 & 42 & 118 & 310 & 779 & 1876 & 4382 & 9948 & 22057 & 47850 & 101844 & \\
\vdots &&&&&&&&&&&&& \end{array}
\]

\end{center}

By taking first the inverse limit of
\[ T_{1,d} \leftarrow T_{2,d} \leftarrow T_{3,d} \leftarrow \, \cdots, \]
one obtains the $\bZ^d$-graded algebra $T_{\infty,d}$. By adapting a
definition of Formanek \cite[p. 52]{EF2}, we refer to $T_{\infty,d}$
as the {\em mixed free trace ring} on $d$ generators.
One can next take the direct limit of
\[ T_{\infty,1} \to T_{\infty,2} \to T_{\infty,3} \to \cdots \]
to obtain the $\bZ^\infty$-graded algebra $T_{\infty,\infty}$,
which is the mixed free trace ring
on countably many generators $x_1,x_2,x_3,\ldots$.

There is a close relationship between Tables 10 and 12
from which we deduce that the following formula apparently holds:
\[ P(C_{\infty,2};t)=(1-2t)P(T_{\infty,2};t). \]
By further heuristic reasoning, one obtains the hypothetical formulae
\[ P(T_{\infty,d};t_1,\ldots,t_d)=\frac{1}{(1-t_1-\cdots-t_d)^2}
\prod_{k\ge2} \frac{1}{1-t_1^k-\cdots-t_d^k} \]
and
\[ P(T_{\infty,\infty};t_1,t_2,\ldots)=\frac{1}{(1-p_1)^2}
\prod_{k\ge2}\frac{1}{1-p_k}. \]

Similarly as in the previous appendix,
the second differences of the coefficients $d_{n,2}(n+k)$
provide yet another sequence:
\[ 1,\, 6,\, 27,\, 103,\, 358,\, 1159,\, \ldots \]

\end{document}